\documentclass{amsart}
\usepackage[ansinew]{inputenc}
\usepackage[english]{babel}
\usepackage{amsmath,latexsym,amssymb,verbatim}

\usepackage[all]{xy}
\usepackage{pdfsync}
\usepackage{hyperref}

\newtheorem{theorem}{Theorem}[section]
\newtheorem{lemma}[theorem]{Lemma}
\newtheorem{proposition}[theorem]{Proposition}
\newtheorem{corollary}[theorem]{Corollary}
\newtheorem{note}[theorem]{Note}
\newtheorem{Formula of adjoint functors}[theorem]{Formula of adjoint functors}

\newtheorem{example}[theorem]{Example}
\newtheorem{examples}[theorem]{Examples}
\newtheorem{definition}[theorem]{Definition}
\newtheorem{notation}[theorem]{Notation}

\DeclareMathOperator{\limi}{{lim}}
\newcommand{\ilim}[1]{\,\underset{#1}{\underset{\to}{\limi}}\,}
\newcommand{\plim}[1]{\,\underset{#1}{\underset{\leftarrow}{\limi}}\,}

\DeclareMathOperator{\Hom}{{Hom}}
\DeclareMathOperator{\Spec}{{Spec}}

\DeclareMathOperator{\Ima}{{Im}}
\DeclareMathOperator{\Der}{{Der}}
\DeclareMathOperator{\charc}{{char}}
\DeclareMathOperator{\Aplic}{{Homeo}}
\DeclareMathOperator{\rad}{{rad}}
\DeclareMathOperator{\Id}{{Id}}
\DeclareMathOperator{\Derinv}{{Derinv}}
\DeclareMathOperator{\tr}{{tr}}

\NoCompileMatrices

\input arrow.tex

\begin{document}

\title{Affine Functors and Duality}

\author{José Navarro}

\author{Carlos Sancho}

\author{Pedro Sancho}


\begin{abstract} A functor of sets $\mathbb X$ over the category of $K$-commutative algebras is said to be an affine functor if its functor of functions, $\mathbb A_{\mathbb X}$, is reflexive and $\mathbb X=\Spec \mathbb A_{\mathbb X}$. We prove that affine functors are equal to a direct limit of affine schemes and that affine schemes, formal schemes, the completion of affine schemes along a closed subscheme, etc., are affine functors.

Endowing an affine functor $\mathbb X$ with a functor of monoids structure 
is equivalent to endowing $\mathbb A_{\mathbb X}$ with  a functor of bialgebras structure.
If $\mathbb G$ is an affine functor of monoids, then $\mathbb A_{\mathbb G}^*$ is the enveloping functor of algebras of $\mathbb G$ and the category of $\mathbb G$-modules is equivalent to the category of $\mathbb A_{\mathbb G}^*$-modules.
Applications of these results include Cartier duality, neutral Tannakian duality for affine group schemes, the equivalence between formal groups and Lie algebras in characteristic zero, etc.
\end{abstract}

\maketitle

\section*{Introduction}

Let $K$ be a (unital associative) commutative ring.
It is well-known that $K$-schemes can be treated as mere ``abstract sets'' by means of their functor of points, which are functors of sets defined over the category of commutative $K$-algebras.
This functorial point of view is particularly useful to study $K$-group schemes and their linear representations (\cite{Demazure2}, \cite{demazure}). Nevertheless, to obtain reasonable results in the case of formal groups, the functors of points have to be endowed with a topology: they are defined over the category of linearly compact $K$-algebras, where $K$ is a pseudo-compact ring (\cite{dieudonne}).

In this paper, that develops ideas introduced in \cite{Amel} and \cite{navarro}, we show that, when considering $K$-modules, linear representations of group schemes, formal groups, etc., as functors over the category of commutative $K$-algebras since the beginning, then the concepts of affine functor and reflexive functor arise naturally, allowing to prove many results from Algebraic Geometry, as obvious consequences of the reflexivity of the functors of modules considered.

 All functors considered in this paper are covariant functors defined over the category of commutative $K$-algebras.

 Given an $K$-module $M$, we denote by
${\mathcal M}$ the functor  ${\mathcal M}(S):= M \otimes_K S$, for all commutative $K$-algebras $S$. We say that $\mathcal M$ is the quasi-coherent module associated with $M$.
If $\mathbb M$ and ${\mathbb M'}$ are functors of $\mathcal K$-modules, then ${\mathbb Hom}_{\mathcal K}(\mathbb M,{\mathbb M'})$
will denote the functor of $\mathcal K$-modules $${\mathbb Hom}_{\mathcal K}(\mathbb M,{\mathbb M'})(S):={\rm
Hom}_{\mathcal S}(\mathbb M_{|S}, {\mathbb M'}_{|S})$$ where $\mathbb M_{|S}$ is the functor $\mathbb M$
restricted to the category of commutative $S$-algebras. We write
$\mathbb M^* := {\mathbb Hom}_{\mathcal K} (\mathbb M ,{\mathcal K})$ and say that this is a dual functor.

The fundamental results on which is based this paper are the reflexivity theorems:

{\bf Theorem (\cite[1.10]{Amel}):} Let $M$ be an $K$-module, then ${\mathcal M}^{**} =
{\mathcal M}$.

{\bf Theorem (\cite[4.4]{navarro}):} Assume that $K$ is a field.
A functor of  $\mathcal K$-modules $\mathbb M$  is reflexive if and only if it is the
inverse limit of its quasi-coherent quotients.

Given a functor of commutative $\mathcal K$-algebras ${\mathbb A}$,
let ${\rm Spec}\, {\mathbb A}$ be the functor: $${\rm Spec}\,
{\mathbb A}(S) := {\rm Hom}_{\mathcal K-alg}({\mathbb A}, {\mathcal S}).$$ In case
$X={\rm Spec}\, A$ is an affine $K$-scheme and $X^\cdot$ stands for its
functor of points ($X^\cdot(S):={\rm Hom}_{K-alg}(A,S)$), then $X^\cdot ={\rm Spec}\, {\mathcal A}$. Let $C$ be a $K$-module,
 if $\mathcal C^*$ is a functor of commutative algebras
 we say that $\Spec\mathcal C^*$ is a formal scheme. If $K$ is a field,  $\Spec\mathcal C^*$ is a formal scheme if and only if  it is a direct limit of finite $K$-schemes (see Note \ref{n3.22}). This is the definition of formal scheme, when $K$ is a field, that can be found in \cite{Demazure2} and \cite{Takeuchi}.

Let $\mathbb X$ be a functor of sets and $\mathbb A_{\mathbb X}:=\mathbb Hom(\mathbb X,\mathcal K)$. We say that $\mathbb X$ is an affine functor if $\mathbb X=\Spec \mathbb A_{\mathbb X}$ and $\mathbb A_{\mathbb X}$ is reflexive.
 We warn the reader that  in the literature affine functors are sometimes defined to be functors of points of affine schemes.

 From now on we will assume, for simplicity, that $K$ is a field. We prove the following theorems.

\begin{theorem} If $\mathbb A$ is a reflexive functor of commutative algebras, $\Spec\mathbb A$ is a direct limit of closed immersions of affine schemes.
\end{theorem}

Let $A$ be a $K$-algebra, we say that $\mathcal A$ is a  quasi-coherent algebra. If
  $\mathbb A$ is the inverse limit of its quasi-coherent algebra quotients we say that
  $\mathbb A$ is a proquasi-coherent algebra.

\begin{theorem} Let $\mathbb A$ be a reflexive functor of commutative algebras. If $\mathbb A$ is a proquasi-coherent algebra, then $\Spec \mathbb A$ is an affine functor and $\mathbb A_{\Spec \mathbb A}=\mathbb A$. If $\mathbb X$ is an affine functor, then $\mathbb A_{\mathbb X}$ is a proquasi-coherent algebra.
\end{theorem}

\begin{theorem}
Affine schemes, formal schemes, the completion of an affine scheme along a closed subscheme are affine functors.
\end{theorem}


\begin{theorem} \label{teo2} An affine functor $\mathbb G$ is a functor of monoids if and only if $\mathbb A_{\mathbb G}$ is a functor of bialgebras, and given two affine functors of monoids $\mathbb G_1$ and $\mathbb G_2$, then
$$\Hom_{mon}(\mathbb G_1,\mathbb G_2)=\Hom_{bialg}(\mathbb A_{\mathbb G_2}, \mathbb A_{\mathbb G_1})$$

 \end{theorem}

In particular, the category of  formal monoids is equivalent to the category of cocommutative bialgebras (see \ref{5.3m}).

In Section \ref{Appendix}, we prove the categorical equivalence of the category of infinitesimal formal groups
 with the category of Lie algebras, and the Poincaré-Birkhoff-Witt Theorem, in characteristic zero (see \cite{Serre}). Let us speak loosely. Let $\mathbb G=\Spec\mathcal A^*$ be an infinitesimal formal group and $\mathcal I^*$ the ideal
of functions of $\mathbb G$ vanishing at the unit element. Then, the natural morphism $S^n(\mathcal I^{*}/\mathcal I^{*2} )\to \mathcal I^{*n}/\mathcal I^{*n+1}$
is an isomorphism, and we prove that the inverse morphism is the morphism induced by the comultiplication morphism $\mathcal A^*\to \mathcal A^*\tilde\otimes\overset n\cdots\tilde\otimes\mathcal A^*$. Dually, if $L=(\mathcal I^{*}/\mathcal I^{*2} )^*=T_e\mathbb G$ is the Lie algebra of $\mathbb G$, $U(L)$ is the universal algebra associated with
$L$ and $U(L)_n:=L\overset n \cdots L$, we obtain that $U(L)_n/U(L)_{n+1}=S^n L$ and as a consequence $A=U(L)$. Given another infinitesimal formal group $\mathbb G'=\Spec\mathcal B^*$, then
$$\aligned {\rm Hom}_{grp}(\mathbb G,{\mathbb G'}) & ={\rm Hom}_{bialg}({\mathcal B^*},{\mathcal A^*})={\rm Hom}_{bialg}
(A,B)\\ & ={\rm Hom}_{bialg}(U(T_e\mathbb G),U(T_e{\mathbb G'}))={\rm Hom}_{Lie}(T_e\mathbb G,T_e{\mathbb G'})\endaligned$$

It is well-known that, for a finite monoid $G$, the category of $K$-linear representations of $G$ is equivalent to the category of $KG$-modules. In \cite[5.4]{Amel} we extended this result to affine group schemes. Now, let $\mathbb G $ be a functor of monoids, such that $\mathbb A_{\mathbb G}$ is reflexive and let $\mathbb G\to \mathbb A_{\mathbb G}^*$ be the natural morphism.
In this paper we prove that the enveloping functor of algebras of $\mathbb G$ is $\mathbb
A_{\mathbb G}^*$, that is,
\begin{equation} \label{envolvente} {\rm Hom}_{mon} (\mathbb G, {\mathbb B}) = {\rm Hom}_{\mathcal K-alg} (\mathbb A_{\mathbb G}^*,{\mathbb B}) \end{equation}
for all dual functors of $\mathcal K$-algebras ${\mathbb
B}$. As consequences of Equality \ref{envolvente}, we obtain the following two theorems.

\begin{theorem}  \label{teo3} The category of dual functors of  $\mathbb G$-modules is equivalent to the category of dual functors of $\mathbb A_{\mathbb G}^*$-modules.
In particular, the category of quasi-coherent ${\mathbb G}$-modules is equi\-valent to the
category of quasi-coherent $\mathbb A_{\mathbb G}^*$-modules.
\end{theorem}

\begin{corollary}
Let
$\mathbb M,\mathbb M'$ be reflexive functors of $\mathbb G$-modules.
 Then, a morphism
    of $\mathcal K$-modules $\mathbb M\to\mathbb M'$ is a morphism of $\mathbb G$-modules if and only if $\mathbb M(K)\to\mathbb M'(K)$ is a morphism of $\mathbb A^*_{\mathbb G}(K)$-modules. Let  $\mathcal M$ be a $\mathbb G$-module, then the set of quasi-coherent $\mathbb G$-submodules of $\mathcal M$ is equal to the set of $\mathbb A_{\mathbb G}^*(K)$-submodules of $M$.

\end{corollary}

In Section \ref{Appendix}, we deduce the Tannaka's characterization of the category of linear representations of an affine group scheme:
Let us talk loosely. If a category of finite-dimensional vector spaces with an extra structure is generated by a unique object $X$, then it is (weak) equivalent to the category of finitely generated $A_X$-modules, for some finite-dimensional algebra  $A_X$. If the category is generated by a set of objects $\{X_i\}$ then it is equivalent to the category of finite generated $\plim{i} \mathcal A_{X_i}$-modules, where $\plim{i} \mathcal A_{X_i}=(\ilim{i} \mathcal A_{X_i}^*)^*=:\mathcal C^*$ is a scheme of algebras.
If in addition, the tensor product of the objects of the category are objects of the category then $\mathcal C^*$ has a comultiplication, that  is, $\mathcal C^*$ is a bialgebra. Therefore, the category is equivalent to the category of finite linear representations of $\Spec C$.

\begin{theorem} \label{teo4} Assume $\mathbb G$ is commutative. Then,
$$\mathbb G^\vee:={\mathbb Hom}_{mon} (\mathbb G, {\mathcal K})=\mathbb Hom_{\mathcal K-alg}(\mathbb A_{\mathbb G}^*,\mathcal K)=\Spec \mathbb A_{\mathbb G}^*
$$
\end{theorem}

As immediate application of Theorem \ref{teo4} and the reflexivity theorem, we deduce the Cartier duality over commutative rings (also see \cite[Ch. I, \S 2, 14]{dieudonne}, where formal schemes are certain functors over the category of commutative linearly compact algebras over a field).

If $G_a=\Spec K[x]$ is the additive group, the category of $G_a^\vee$ modules is equivalent
to the category of $K[x]$-modules. Then,
$$\Hom_{mon}(G_a^\vee, {\mathbb End}_{K}(V))=\mathbb End_K(V)$$
If $V$ is a $K$-algebra, we prove that
$$\Hom_{mon}(G_a^\vee, {\mathbb End}_{K-alg}(V))=\Der_K(V,V)$$
More generally, we prove  (\ref{8.4}) that if $G$ is $K$-group
scheme, then
$$ \Hom_{mon}(G_a^\vee,G)=T_{e}G$$
and if  $X$ is a $K$-scheme,  then
$\Hom_{mon}(G_a^\vee,{\mathbb End}\, X)=\Der X$ (\ref{55}).

That is, giving a vector field $D$ on a $K$-scheme $X$ is equivalent to giving a morphism
of functors of groups $exp_D\colon G_a^\vee\to {\mathbb End}\, X$; and giving an invariant vector field on a $K$-group scheme $G$ is equivalent to giving a morphism
of functors of groups $exp_D\colon G_a^\vee\to G$.

Every morphism $\Spec \mathcal C^*\to X$, from a infinitesimal formal scheme (i.e., $C^*$ is a local algebra) to a $K$-scheme uniquely factors via $\Spec C^*$, that is,

$$\Hom(\Spec {\mathcal C}^*, X)=\Hom_{K-sch}(\Spec {C}^*, X).$$
Then, $exp_D$ uniquely factors through $\Spec K[x]^*$. In characteristic zero, $K[x]^*\simeq K[[z]]$ and one has the classical exponential map $exp_D\colon \Spec K[[z]]\to G$,  associated with $D$ (\cite[II  6 3]{Demazure}). In characteristic $p>0$,
one has an isomorphism of $K$-schemes $\Spec K[x]^*=\Spec K[[x_0,\ldots,x_n,\ldots]]/(x_0^p,\ldots,x_n^p,\ldots)$
(see \ref{7.9}). We apply this construction to prove the existence and uniqueness
of ``analytic'' solutions of an algebraic differential equation in
arbitrary characteristic (see \ref{632}):

Let $\delta$ be the invariant field on $\Spec K[x]^*$. That is, in characteristic zero $\delta=\partial_z$, in characteristic $p>0$,
$\delta=\partial_{x_0}+x_0^{p-1}/(p-1)!\partial_{x_1}
+x_0^{p-1}x_1^{p-1}/(p-1)!^2\partial_{x_2}+\cdots$.

\begin{theorem} Let $X$ be a $K$-scheme, let $y \in X$ be a rational point and let $D$ be a vector field on $X$. Then,  $exp_{D,y}\colon  \Spec K[x]^* \to X$, $exp_{D,y}(x):=exp_{D}(x)$ is
the only morphism $f\colon \Spec K[x]^* \to X$ such that $f(0)=y$ and $f(\delta_{\mu})=D_{f(\mu)}$, for every point
$\mu\in (\Spec K[x]^*)(S)$.\end{theorem}

If $X$ is a complete algebraic variety, then the scheme-theoretic image of $exp_D$ is a commutative algebraic group,
that we define to be the algebraic group associated with $D$ (\ref{75}). The minimal subvariety tangent to $D$ passing
through a point is the orbit of the point under the action of the algebraic group associated with $D$.

Let $\Spec \mathcal C^*$ be a formal monoid and $D_{\mathcal C}=\{w\in C^*\colon w(I)=0$ for some bilateral ideal $I\subset C$ of finite codimension$\}\subset C^*$.
Then, $\Spec D_{\mathcal C}$ is an affine monoid scheme,
because
$\mathcal D_{\mathcal C}^*=\plim{I}\mathcal C/\mathcal I$ is a scheme of bialgebras, and
\begin{equation} \label{i} \aligned \Hom_{mon}(\Spec \mathcal C^*, \Spec A)& =\Hom_{bialg}(\mathcal A,\mathcal C^*)=\Hom_{bialg}(\mathcal C, {\mathcal A}^*)\\ & =
\Hom_{bialg}(\plim{I}\mathcal C/\mathcal I,\mathcal A^*)=\Hom_{bialg}({\mathcal A},\mathcal D_{\mathcal C} )\\ & =
\Hom_{mon}( \Spec D_{\mathcal C},\Spec A)\endaligned\end{equation}

In the case the algebraic group associated with $D$ is an affine algebraic
group,  then it is isomorphic to a quotient of $\Spec D_{\mathcal K[x]}$, by Equation \ref{i}.
Then, by Theorem \ref{37}, it is isomorphic to
$$\begin{array}{ll} G_a^\delta\times G_m^n,\,
\delta=0,1,\qquad &\text{ if }\charc k=0 \\ \alpha_r^\vee\times \mu_1^n, \qquad \ \, & \text{ if } \charc k=p>0\end{array}$$

Finally,  in Section \ref{76}, we calculate the algebraic group associated with a field on $\mathbb P ^n(k)$, where $k$ is an
algebraically closed field of arbitrary characteristic, recovering in this way the results from \cite{Juan} for the
case of characteristic zero.

\begin{theorem} Let $K$ be an algebraically closed field and let $\pi\colon \mathbb A^{n}(K) \backslash 0\to \mathbb P^{n-1}(K)$ be the projectivization map. Let
$D=\pi(\sum_{ij}\lambda_{ij}x_i\partial_{x_j})$ be a vector field on $\mathbb P^{n-1}(K)$ and let $G$ be its associated
algebraic group.

\begin{enumerate}

\item If $\charc K=0$
$$G\simeq G_m^r\times G_a^\delta$$
where $r$ is the dimension of the $\mathbb Q$-affine space
generated by the eigenvalues of the matrix  $(\lambda_{ij})$ in
$K$, $\delta=0$ in case the matrix is diagonalizable and
$\delta=1$ otherwise.

\item If $\charc K=p>0$
$$G\simeq \mu_1^r\times \alpha_{m+1}^\vee$$
where $r$ is the dimension of the $\mathbb Z/p\mathbb Z$-affine space generated by the eigenvalues of the matrix $(\lambda_{ij})$
in $K$, and $m$ is such that, if $s$ is the greatest of the orders of the Jordan boxes, then $p^m\leq s-1<p^{m+1}$ (if
$s=1$ we say that $m=-1$).\end{enumerate}\end{theorem}

\section{Main notations}

$R$ is a (unital associative) commutative ring, $S$ is a commutative $R$-algebra.
$A,B$ are $R$-algebras (or $R$-bialgebras).
$M,N$ are $R$-modules.
$K$ is a field, $V$ is a $K$-vector space or a free $R$-module.

$\mathbb X, \mathbb Y$ are (covariant) functors of sets. $\mathbb G$ is a functor of monoids (or semigroups, or groups).
$\mathbb M$, $\mathcal M,\mathcal N$ are functors of $\mathcal R$-modules.
$\mathcal M,\mathcal N$ are called quasi-coherent modules, $\mathcal M(S):=M\otimes_RS$.

$\Hom_{\mathcal R}(\mathbb M, \mathbb M')$ is the set of morphisms of functors of $\mathcal R$-modules from $\mathbb M$ to $\mathbb M'$. $\mathbb Hom_{\mathcal R}(\mathbb M, \mathbb N)$ is the $\mathcal R$-module functor of morphisms of $\mathcal R$-modules.
$\mathbb M^*:= \mathbb Hom_{\mathcal R}(\mathbb M, \mathcal R)$ is the dual functor (of $\mathbb M$). The dual functor of $\mathcal M$, $\mathcal M^*$ is called module scheme.

$\mathbb A,\mathbb B$, $\mathcal A,\mathcal B$, $\mathcal C^*$ are functors of $\mathcal R$-algebras or $\mathcal R$-bialgebras ($\mathcal A,\mathcal B$ are quasi-coherent $\mathcal R$-modules and
$\mathcal C^*$ is an $\mathcal R$-module scheme).
$\Spec A=\Spec \mathcal A$ is an affine scheme.
$\Spec \mathcal C^*$  is called formal scheme.

$\mathbb A_{\mathbb X}:=\mathbb Hom(\mathbb X,\mathcal R)$ is the functor of functions of the functor of sets $\mathbb X$.

$\mathfrak F$ is a wide family of reflexive $\mathcal R$-modules, which contains free quasi-coherent $\mathcal R$-modules and it is closed by the functor
$\mathbb Hom_{\mathcal R}(-,-)$ and ``essentially'' closed by the functor $-\otimes_{\mathcal R}-$.

\section{Preliminaries}

Let $R$ be a commutative ring (associative with a unit). All functors considered in this paper are covariant functors over the category of commutative $R$-algebras (always assumed to be associative with a unit). A functor $\mathbb X$ is said to be a functor of sets (resp. monoids, etc.) if $\mathbb X$ is a functor from the category of  commutative $R$-algebras to the category of sets (resp. monoids, etc.).

\begin{notation} For simplicity, given a functor of sets $\mathbb X$,
we sometimes use $x \in \mathbb X$  to denote $x \in \mathbb X(S)$. Given $x \in \mathbb X(S)$ and a morphism of commutative $R$-algebras $S \to S'$, we still denote by $x$ its image by the morphism $\mathbb X(S) \to \mathbb X(S')$.\end{notation}

Let $\mathcal R$ be the functor of rings defined by ${\mathcal R}(S):=S$, for all commutative $R$-algebras $S$.
A functor of sets $\mathbb M$ is said to be a functor of $\mathcal R$-modules if we have morphisms of functors of sets, $\mathbb M\times \mathbb M\to \mathbb M$ and ${\mathcal R}\times \mathbb M\to \mathbb M$, so that
$\mathbb M(S)$ is an $S$-module, for every commutative $R$-algebra $S$.
A functor of algebras (associative with a unit), $\mathbb A$,
is said to be a functor of $\mathcal R$-algebras if we have a morphism of functors of algebras
$\mathcal R\to \mathbb A$
(and $\mathcal R(S)=S$ commutes with all the elements of $\mathbb A(S)$, for every commutative $R$-algebra $S$).

Given a commutative $R$-algebra $S$, we denote by $\mathbb M_{|S}$ the functor $\mathbb M$ restricted to the category of commutative $S$-algebras.

 Let $\mathbb M$ and $\mathbb M'$ be functors of $\mathcal R$-modules.
 A morphism of functors of $\mathcal R$-modules $f\colon \mathbb M\to \mathbb M'$
 is a morphism of functors such that the defined morphisms $f_S\colon \mathbb M(S)\to
 \mathbb M'(S)$ are morphisms of $S$-modules, for all commutative $R$-algebras $S$.
 We will denote by $\Hom_{\mathcal R}(\mathbb M,\mathbb M')$ the  set of all the morphisms of $\mathcal R$-modules from $\mathbb M$ to $\mathbb M'$.
We will denote by ${\mathbb Hom}_{\mathcal R}(\mathbb M,\mathbb M')$\footnote{In this paper, we will only  consider functors $\mathbb M$ and $\mathbb M'$ such that $\Hom_{\mathcal S}(\mathbb M_{|S},\mathbb {M'}_{|S})$ are sets, for all $S$.} the functor of $\mathcal R$-modules $${\mathbb Hom}_{\mathcal R}(\mathbb M,\mathbb M')(S):={\rm Hom}_{\mathcal S}(\mathbb M_{|S}, \mathbb M'_{|S})$$  Obviously,
$$(\mathbb Hom_{\mathcal R}(\mathbb M,\mathbb M'))_{|S}=
\mathbb Hom_{\mathcal S}(\mathbb M_{|S},\mathbb M'_{|S})$$

\begin{notation} We denote $\mathbb M^*=\mathbb Hom_{\mathcal R}(\mathbb M,\mathcal R)$.\end{notation}

\begin{notation} Tensor products, direct limits, inverse limits, kernels, cokernels, images, etc., of functors of $\mathcal R$-modules are regarded in the category of functors of $\mathcal R$-modules.\end{notation}

\begin{definition} Given an $R$-module $M$, the functor of $\mathcal R$-modules ${\mathcal M}$  defined by ${\mathcal M}(S) := M \otimes_R S$ is called a quasi-coherent $\mathcal R$-module.
\end{definition}

\begin{proposition} \cite[1.3]{Amel}\label{tercer}
For every functor of ${\mathcal R}$-modules $\mathbb M$ and every $R$-module $M$, it
holds that
$${\rm Hom}_{\mathcal R} ({\mathcal M}, \mathbb M) = {\rm Hom}_R (M, \mathbb M(R))$$
\end{proposition}

The functors $M \rightsquigarrow {\mathcal M}$, ${\mathcal M} \rightsquigarrow {\mathcal M}(R)=M$ establish an equivalence between the category of $\mathcal R$-modules and the category of quasi-coherent $\mathcal R$-modules (\cite[1.12]{Amel}). In particular, ${\rm Hom}_{\mathcal R} ({\mathcal M},{\mathcal M'}) = {\rm Hom}_R (M,M')$.
 For any pair of $R$-modules $M$ and $N$, the quasi-coherent module associated with $M\otimes_R N$ is $\mathcal M\otimes_{\mathcal R}\mathcal N$. ${\mathcal M}_{\mid S}$ is the quasi-coherent $\mathcal S$-module associated with $M \otimes_R
S$

The functor ${\mathcal M}^* = {\mathbb Hom}_{\mathcal R} ({\mathcal M}, {\mathcal R})$ is  called an $\mathcal R$-module scheme. Moreover, ${\mathcal M}^*(S)={\rm Hom}_S(M\otimes_RS,S)={\rm Hom}_R(M,S)$ and it is easy to check that $(\mathcal M^*)_{|S}$ is an $\mathcal S$-module scheme.

\begin{proposition}  \cite[1.8]{Amel} \label{a1.8}
Let $M$, $M'$ be $R$-modules. Then $${\mathbb Hom}_{\mathcal R} ({\mathcal M^*}, {\mathcal M'}) = {\mathcal M} \otimes_{\mathcal R} {\mathcal M'}$$
\end{proposition}

As a corollary we obtain the following theorem.

\begin{theorem} \cite[1.10]{Amel}\label{reflex}
Let $M$ be an $R$-module. Then $${\mathcal M^{**}} = {\mathcal M}$$
\end{theorem}

 The functors $\mathcal M\rightsquigarrow \mathcal M^*$ and $\mathcal M^{*}\rightsquigarrow \mathcal M^{**}=\mathcal M$ establish an anti-equivalence between the categories of quasi-coherent modules
and module schemes.

Let us recall the Formula of adjoint functors.

\begin{definition} Let $i^*\colon R\to S$ be a commutative $R$-algebra.
Given a functor of $\mathcal R$-modules, $\mathbb M$, let $i^* \mathbb M$  be the functor of
$\mathcal S$-modules defined by
$(i^* \mathbb M)(S') := \mathbb M(S')$.
Given a functor of $\mathcal S$-modules, $\mathbb N$, let $i_* \mathbb N$ be the functor of
$\mathcal R$-modules defined by
$(i_* \mathbb N)(R') := \mathbb N(S \otimes_R R')$.
\end{definition}

\begin{Formula of adjoint functors} \cite[2.11]{navarro}\label{adj}
Let $\mathbb M$ be a functor of ${\mathcal R}$-modules and let $\mathbb N$ be a functor of  $\mathcal S$-modules. Then, it holds that
$${\rm Hom}_{\mathcal S} (i^*\mathbb M, \mathbb N) = {\rm Hom}_{\mathcal R} (\mathbb M, i_*\mathbb N)$$
\end{Formula of adjoint functors}

\begin{corollary} \cite[2.12]{navarro} \label{adj2} Let $\mathbb M$ be a functor of $\mathcal R$-modules. Then
$$\mathbb M^*(S)=\Hom_{\mathcal R}(\mathbb M,\mathcal S)$$
for all commutative $R$-algebras $S$.\end{corollary}

\begin{definition} Let $\mathbb M$ be a functor of $\mathcal R$-modules. We will say that
$\mathbb M^*$ is a dual functor.
We will say that a functor of $\mathcal R$-modules ${\mathbb M}$ is reflexive if ${\mathbb M}={\mathbb M}^{**}$.\end{definition}

\begin{examples}  Quasi-coherent modules and module schemes are reflexive functors of $\mathcal R$-modules.\end{examples}

\begin{proposition} \cite[2.16]{navarro} \label{n2.16}
Let $\mathbb A$ be a functor of $\mathcal R$-algebras such that $\mathbb A^*$ is a
reflexive functor of $\mathcal R$-modules. The closure of dual functors of
$\mathcal R$-algebras of $\mathbb A$ is $\mathbb A^{**}$, that is, it holds the functorial
equality $${\rm Hom}_{\mathcal R-alg}(\mathbb A,\mathbb B)={\rm Hom}_{\mathcal R-alg}(\mathbb A^{**},\mathbb B)$$
for every dual functor of $\mathcal R$-algebras $\mathbb B$.

Moreover, endowing a dual functor of $\mathcal R$-modules $\mathbb M^*$ with a structure of ${\mathbb A}$-module is
equivalent to endowing $\mathbb M^*$ with a structure of $\mathbb
A^{**}$-module.
\end{proposition}

 \begin{definition} \cite[5.2]{navarro} \label{n5.2} Let $\mathfrak F$ be the family of dual functors of $\mathcal R$-modules, $\mathbb M$, such that there exist  a set $J$ (which depends on $\mathbb M$), a structure of functor of $\prod_J \mathcal R$-modules on $\mathbb M$  and inclusions of functors of $\prod_J \mathcal R$-modules
$$\oplus_J\mathcal R\subseteq \mathbb M\subseteq \prod_J \mathcal R$$
\end{definition}

\begin{proposition}  \cite[5.3,5.8,5.9]{navarro} \label{n5.9} Every $\mathbb M\in \mathfrak F$ is a functor of $\mathcal R$-modules reflexive.
If $M$ is a free $R$-module then $\mathcal M,\mathcal M^*\in\mathfrak F$. If $\mathbb M,\mathbb M'\in\mathfrak F$, then $\mathbb Hom_{\mathcal R}(\mathbb M,\mathbb M')\in\mathfrak F$ and
$(\mathbb M\otimes_{\mathcal R}\mathbb M')^{**}\in\mathfrak F$, which satisfies
     $$
     \Hom_{\mathcal R}((\mathbb M\otimes_{\mathcal R}\mathbb M')^{**},\mathbb {M''})=\Hom_{\mathcal R}(\mathbb M\otimes_{\mathcal R}\mathbb M',\mathbb {M''})$$
for every reflexive  functor of $\mathcal R$-modules, $\mathbb {M''}$.

\end{proposition}

\begin{proposition}  \cite[Section 1]{navarro} \label{nS1} Let $\mathbb A,\mathbb M,\mathbb M'$ be reflexive functors of $\mathcal R$-modules.
Assume that $\mathbb A,\mathbb M,\mathbb M' \in \mathfrak F$ or  that $R=K$ is a field.
If $\mathbb A$ is a functor of $\mathcal R$-algebras and
    $\mathbb M,\mathbb M'$ are functors of $\mathbb A$-modules, then a morphism
    of $\mathcal R$-modules $\mathbb M\to\mathbb M'$ is a morphism of $\mathbb A$-modules if and only if $\mathbb M(R)\to\mathbb M'(R)$ is a morphism of $\mathbb A(R)$-modules. Let  $\mathcal M$ be an $\mathbb A$-module, then the set of quasi-coherent $\mathbb A$-submodules of $\mathcal M$ is equal to the set of $\mathbb A(R)$-submodules of $M$.
\end{proposition}

\begin{definition} Let $A$ be a $R$-algebra, we say that $\mathcal A$ is a  quasi-coherent algebra. If
  $\mathbb A$ is the inverse limit of its quasi-coherent algebra quotients we say that
  $\mathbb A$ is a proquasi-coherent algebra.
\end{definition}

\begin{proposition} \cite[3.18, 5.17]{navarro} \label{n5.17}
Let $\mathbb A$ be a reflexive functor of $\mathcal R$-algebras. Assume that $\mathbb A\in\mathfrak F$ or that $R=K$ is a field.  Every morphism of $\mathcal R$-algebras $\phi\colon \mathbb A\to \mathcal B$ uniquely factors  through an epimorphism of functors of algebras onto the quasi-coherent algebra associated with $\Ima \phi_R$.
Then, if $\{\mathcal A_i\}_i$ is the set of the quasi-coherent algebra quotients of $\mathbb A$,
$$\Hom_{\mathcal R-alg}(\mathbb A,\mathbb B)=\ilim{i} \Hom_{\mathcal R-alg}(\mathcal A_i,\mathbb B)$$
for every functor of proquasi-coherent algebras $\mathbb B$.

\end{proposition}

\begin{note} \label{n3.2} Let $\mathcal C^*\in\mathfrak F$ be a functor of $\mathcal R$-algebras.
$\mathcal C^*$ is a proquasi-coherent algebra and the quasi-coherent algebra quotients of
$\mathcal C^*$ are $R$-modules of finite type (see proof of Proposition  \cite[5.20]{navarro}).
\end{note}

\begin{notation} Given two $\mathcal R$-modules, $\mathbb M$ and $\mathbb M'$, we denote
$\mathbb M\tilde\otimes \mathbb M':=(\mathbb M^*\otimes\mathbb M'^*)^*$.
\end{notation}

\begin{proposition} \label{2.21} Let $\mathbb M,\mathbb M'\in\mathfrak F$. By Proposition \ref{n5.9},
$\mathbb M\tilde\otimes \mathbb M'=(\mathbb M^*\otimes\mathbb M'^*)^*=\mathbb Hom_{\mathcal R}(\mathbb M^*,\mathbb M')\in\mathfrak F$. Given two modules $M$ and $N$,
$$\mathcal M\tilde \otimes\mathcal N=(\mathcal M^*\otimes\mathcal N^*)^*=
\mathbb Hom_{\mathcal R}(\mathcal M^*,\mathcal N)\overset{\text{\ref{a1.8}}}=\mathcal M\otimes\mathcal N$$
Finally, $\mathcal M^*\tilde\otimes\mathcal N^*=(\mathcal M\otimes \mathcal N)^*$ is a module scheme and
$$\Hom_{\mathcal R}(\mathcal M^*\tilde\otimes\mathcal N^*,\mathbb P^*)=
\Hom_{\mathcal R}(\mathbb P, \mathcal M\otimes\mathcal N)=
\Hom_{\mathcal R}(\mathcal M^*\otimes\mathcal N^*,\mathbb P^*)$$
for all dual modules $\mathbb P^*$.

\end{proposition}

\begin{proposition} \cite[5.22]{navarro} \label{n5.22}
Let $\mathbb A, \mathbb B\in \mathfrak F$ be functors of proquasi-coherent algebras. Then,
$\mathbb A\tilde\otimes \mathbb B:=(\mathbb A^*\otimes \mathbb B^*)^*\in\mathfrak F$ is a proquasi-coherent algebra such that
$$\Hom_{\mathcal R-alg}(\mathbb A\otimes \mathbb B,\mathbb C)=\Hom_{\mathcal R-alg}((\mathbb A^*\otimes \mathbb B^*)^*,\mathbb C)$$
for every functor of proquasi-coherent algebras $\mathbb C$.
\end{proposition}

\begin{definition}  A functor  of proquasi-coherent algebras $\mathbb B\in \mathfrak F$
is said to be a functor of  bialgebras  (resp. a functor of proquasi-coherent bialgebras) if $\mathbb B^*$ is a functor of $\mathcal R$-algebras (resp. a functor of proquasi-coherent $\mathcal R$-algebras) such that the dual morphisms of the multiplication morphism $m\colon \mathbb B^*\otimes \mathbb B^*\to \mathbb B^*$ and the unit morphism $u\colon \mathcal R\to \mathbb B^*$ are morphisms of functors of $\mathcal R$-algebras.

Let $\mathbb B,\mathbb B'$ be two functors of bialgebras. We will say that a morphism of $\mathcal R$-modules, $f\colon \mathbb B\to\mathbb B'$ is a morphism of functors of bialgebras if $f$ and $f^* \colon {\mathbb B'}^*\to\mathbb B^*$ are morphisms of functors of $\mathcal R$-algebras.

\end{definition}

\begin{theorem} \cite[5.27]{navarro} \label{n5.27} Let ${\mathcal C}_{\mathfrak F-Bialg.}$ be the category of
functors $\mathbb B\in\mathfrak F$ of proquasi-coherent bialgebras. The functor ${\mathcal C}_{\mathfrak F-Bialg.}\rightsquigarrow{\mathcal C}_{\mathfrak F-Bialg.}$, $\mathbb B \rightsquigarrow {\mathbb B}^*$ is a categorical anti-equivalence.
\end{theorem}

\begin{notation} Let $\mathbb A$ be a reflexive functor of $\mathcal K$-algebras and
let $\{\mathcal A_i\}$ the set of quasi-coherent quotients of $\mathbb A$ such that
$\dim_K A_i<\infty$. We denote $\bar{\mathbb A}:=\plim{i} \mathcal A_i$ which
is an algebra scheme because $\mathcal A_i^*$ is quasi-coherent and $\plim{i} \mathcal A_i=(\ilim{i} \mathcal A_i^*)^*$.\end{notation}

\begin{proposition} \cite[5.9]{Amel} Let $\mathbb A$ be a reflexive functor of $\mathcal K$-algebras. Then,
$$\Hom_{\mathcal K-alg}(\mathbb A,\mathcal C^*)=
\Hom_{\mathcal K-alg}(\bar{\mathbb A},\mathcal C^*)$$
for all algebra schemes $\mathcal C^*$.
\end{proposition}

\begin{theorem} \cite[5.30]{navarro} \label{n5.30} Let $\mathbb B\in\mathfrak F$ be a functor of proquasi-coherent $\mathcal K$-bialgebras. Then, $\bar{\mathbb B}$ is a scheme of bialgebras and

$$\Hom_{bialg}(\mathbb B,\mathcal C^*)=
\Hom_{bialg}(\bar{\mathbb B},\mathcal C^*)$$
for all bialgebra schemes $\mathcal C^*$.
\end{theorem}

\section{Affine functors}

Let $X$ be an $R$-scheme and let $X^\cdot $ be
the functor of points of $X$; i.e., $X^\cdot$ is the functor of sets  $$X^\cdot (S) ={\rm
Hom}_{R-{\rm sch}}({\rm Spec}\, S,X)$$  For any other scheme $Y $, Yoneda's
lemma proves that $${\rm Hom}_{R-sch}(X,Y)={\rm Hom} (X^\cdot ,Y^\cdot ),$$ so $X^\cdot  \simeq Y^\cdot $ if and only if $X \simeq Y$. We will sometimes denote $X^\cdot  =X$. If $X={\rm Spec}\, A$ is an affine scheme, then
$$({\rm Spec}\, A)^\cdot (S) ={\rm Hom}_{R-{\rm sch}}({\rm Spec}\, S,{\rm Spec}\, A)={\rm Hom}_{R-alg} (A,S)$$

\begin{definition}
Given a functor of commutative $\mathcal R$-algebras $\mathbb {A}$, the functor
${\rm Spec}\, \mathbb {A}$, ``spectrum of $\mathbb {A}$'', is defined
to be $$({\rm Spec}\, \mathbb {A})(S) := {\rm Hom}_{\mathcal R-alg}
(\mathbb {A}, {\mathcal S})$$ for every commutative $R$-algebra $S$.
\end{definition}

\begin{proposition}
Let $\mathbb A$ be a functor of commutative algebras. Then, $${\rm
Spec}\,{\mathbb A}={\mathbb Hom}_{\mathcal R-alg}({\mathbb A},{\mathcal R}).$$
\end{proposition}

\begin{proof}
By Adjoint Formula (\cite[2.11]{navarro}), restricted to
the morphisms of algebras, it holds that
$${\mathbb Hom}_{\mathcal R-alg} ({\mathbb A}, {\mathcal R}) (S)
= {\rm Hom}_{\mathcal S-alg} ({\mathbb A}_{|S}, {\mathcal S})
= {\rm Hom}_{\mathcal R-alg} ({\mathbb A}, {\mathcal S})  = ({\rm Spec}\, {\mathbb A}) (S).$$
\end{proof}

Therefore, ${\rm Spec}\,{\mathbb A}={\mathbb Hom}_{\mathcal R-alg} ({\mathbb A},{\mathcal
R}) \subset {\mathbb Hom}_{\mathcal R}({\mathbb A},{\mathcal R}) = {\mathbb A}^*$.

\begin{notation} Given a functor of sets $\mathbb X$, the functor $\mathbb A_{\mathbb X}:={\mathbb Hom}({\mathbb X}, {\mathcal
R})$ is said to be the functor  of functions of $\mathbb X$.\end{notation}

\begin{proposition} \label{homspe} Let $\mathbb X$ be a functor of sets and $\mathbb A_{\mathbb X}$ its functor of functions. Then,

$${\rm Hom}(\mathbb X,\Spec \mathbb B)={\rm Hom}_{\mathcal R-alg}(\mathbb B,\mathbb A_{\mathbb X}),$$ for every functor of commutative algebras, $\mathbb B$.

\end{proposition}

\begin{proof} Given $f\colon \mathbb X\to \Spec\mathbb B$, let  $f^*\colon \mathbb B\to \mathbb A_{\mathbb X}$ be defined by $f^*(b)(x):=f(x)(b)$, for every $x\in \mathbb X$. Given $\phi\colon \mathbb B\to \mathbb A_{\mathbb X}$, let
$\phi^*\colon \mathbb X\to \Spec\mathbb B$ be defined by $\phi^*(x)(b):=
\phi(b)(x)$, for all $b\in\mathbb B$. It is easy to check that $f=f^{**}$ and $\phi=\phi^{**}$.

\end{proof}

\begin{example} \label{ejemplo}
If $A$ is a commutative $R$-algebra, then ${\rm Spec}\, {\mathcal A} =
({\rm Spec}\, A)^\cdot$ and $\mathbb A_{\Spec\mathcal A}={\mathbb Hom}(\Spec\mathcal A,\mathcal R)={\mathbb Hom}(\Spec\mathcal A,\Spec R[x])={\mathbb Hom}_{\mathcal R-alg}(\mathcal R[x],\mathcal A)=\mathcal A$.
\end{example}

If $R=K$ is a field and $X$ is a noetherian $K$-scheme, then the functor of functions of $X^\cdot$ is a quasi-coherent $\mathcal R$-module.

\begin{definition} We will say that a functor of sets $\mathbb X$ is affine
when $\mathbb X=\Spec\mathbb A_{\mathbb X}$ and $\mathbb A_{\mathbb X}$ is reflexive.
\end{definition}

Let $\mathbb X$ and $\mathbb Y$ be affine functors. By Proposition \ref{homspe},
$$\Hom(\mathbb X,\mathbb Y)=\Hom_{\mathcal R-alg}(\mathbb A_{\mathbb Y},\mathbb A_{\mathbb X})$$

\begin{example} \label{eje4.5} Affine schemes, $\Spec \mathcal A$, are affine functors, by Example
\ref{ejemplo}.\end{example}
\begin{proposition} \label{3.10} Let $\mathbb X = \ilim{i}\Spec\mathcal A_i$.  Then,
$\mathbb A_{\mathbb X}=\plim{i}\mathcal A_i$.
\end{proposition}

\begin{proof}
It holds that  ${\rm Hom} (({\rm Spec}\, A)^\cdot, {\mathbb Y}) =
{\mathbb Y}(A)$ for every functor of sets ${\mathbb Y}$, by Yoneda's lemma.
Then,

$$\mathbb A_{\mathbb X}={\mathbb Hom}(\mathbb X,\mathcal R)   ={\mathbb  Hom}
(\underset{\underset{i}{\longrightarrow}} {\lim}\, {\rm Spec}\, {\mathcal
A}_i, \mathcal R)  =
\underset{\underset{i}{\longleftarrow}}{\lim}\,{\mathbb  Hom}
({\rm Spec}\, {\mathcal A}_i,\mathcal R)
= \plim{i}\mathcal A_i$$
\end{proof}

\begin{theorem} \label{t1.9} Let $\mathbb A$ be a reflexive functor of commutative algebras. Assume $R=K$ is a field or $\mathbb A\in\mathfrak F$. Then,
$\Spec \mathbb A$ is a direct limit of closed immersions of affine schemes.\end{theorem}

\begin{proof} Let $\{\mathcal A_i\}_i$ be the set of of commutative quasi-coherent algebra quotients of $\mathbb A$.
$\Spec\mathbb A=\mathbb Hom_{\mathcal R-alg}(\mathbb A,\mathcal R)\overset{\text{\ref{n5.17}}}=
\ilim{i} \mathbb Hom_{\mathcal R-alg}(\mathcal A_i,\mathcal R)=\ilim{i} \Spec \mathcal A_i$.\end{proof}

\begin{theorem} \label{t1.92} Let $\mathbb A$ be a reflexive functor of commutative algebras. Assume $R=K$ is a field or $\mathbb A\in\mathfrak F$.
If $\mathbb A$ is a proquasi-coherent algebra then
$\Spec\mathbb A$ is an affine functor and $\mathbb A=\mathbb A_{\Spec\mathbb A}$.
If $\Spec \mathbb A$ is an affine functor then
$\mathbb A_{\Spec\mathbb A}$ is a proquasi-coherent algebra.
\end{theorem}

\begin{proof} Let $\{\mathcal A_i\}_i$ be the set of of commutative quasi-coherent algebra quotients of $\mathbb A$.
Then, $\mathbb A_{\Spec\mathbb A}=\plim{i}\mathcal A_i$, by Theorem \ref{t1.9} and Proposition \ref{3.10}.

If $\mathbb A$ is a proquasi-coherent algebra, then $\mathbb A=\plim{i}\mathcal A_i=\mathbb A_{\Spec\mathbb A}$, and $\Spec\mathbb A$ is an affine functor.

Suppose that $\Spec \mathbb A$ is an affine functor.
Let $f\colon \mathbb A_{\Spec\mathbb A} \to \mathcal B$ be a morphism of functors of algebras. The composition morphism $\mathbb A\to \mathbb A_{\Spec\mathbb A} \to \mathcal B$ factors through some $\mathcal A_i$. As $\Spec \mathbb A=\Spec \mathbb A_{\Spec\mathbb A}$, $f$ factors through $\mathcal A_i$. In conclusion, the set of quasi-coherent algebra quotients of $\mathbb A_{\Spec\mathbb A}$ is $\{\mathcal A_i\}_{i\in I}$, and $\mathbb A_{\Spec\mathbb A}$ is a proquasi-coherent algebra.

\end{proof}

\begin{definition} Let $\mathbb X$ be a functor of sets.
 Let $\mathcal R\mathbb X$ be the functor of $\mathcal R$-modules defined by $$\mathcal R\mathbb X(S):=\oplus_{\mathbb X(S)} S=\{ \mbox{formal finite $S$-linear  combinations of elements of } \mathbb X(S) \}$$\end{definition}

 Clearly, ${\mathbb Hom}(\mathbb X,{\mathbb M})=
{\mathbb Hom}_{\mathcal R}(\mathcal R\mathbb X,{\mathbb M})$, for all functors of $\mathcal R$-modules, $\mathbb M$.

Observe that $\mathbb A_{\mathbb X}={\mathbb Hom}(\mathbb X,{\mathcal R})=(\mathcal R\mathbb X)^*$ is a dual functor.

\begin{proposition}\label{2.3} Let $\mathbb X$ be a functor of sets. Let $\mathbb B_{\mathbb X}$ be a
functor of $\mathcal R$-modules such that $\mathbb A_{\mathbb X}=\mathbb B_{\mathbb X}^*$. Then,
$${\rm Hom}({\mathbb X}, {\mathbb M^*}) =  {\rm Hom}_{\mathcal R}(\mathbb B_{\mathbb X}, {\mathbb M^*})$$ for every dual functor of $\mathcal R$-modules ${\mathbb M}^*$.
\end{proposition}

\begin{proof} It holds that
 $$\aligned {\rm Hom}({\mathbb X}, {\mathbb M}^*) & =  {\rm Hom}_{\mathcal R}(\mathcal R\mathbb X, {\mathbb M}^*) ={\rm Hom}_{\mathcal R}(\mathcal R\mathbb X\otimes{\mathbb M},\mathcal R) ={\rm Hom}_{\mathcal R}({\mathbb M}, {\mathbb A}_{\mathbb X})\\ & ={\rm Hom}_{\mathcal R}(\mathbb B_{\mathbb X}, {\mathbb M}^*)\endaligned$$

\end{proof}

\begin{proposition} \label{244} Let $\mathbb X_i$, $i=1,\ldots,n$  be functors of sets and $\mathbb A_{\mathbb X_i}=\mathbb B_{\mathbb X_i}^*$. Then,
$$\mathbb Hom(\mathbb X_1\times \cdots\times \mathbb X_n,\mathbb M^*)=\mathbb Hom_{\mathcal R}(\mathbb B_{\mathbb X_1}\otimes \cdots\otimes \mathbb B_{\mathbb X_n},\mathbb M^*)$$
for all dual functors of $\mathcal R$-modules $\mathbb M^*$. In particular, if $\mathbb X_i$ are affine functors, then

$$\mathbb A_{\mathbb X_1\times \cdots\times\mathbb X_n}=(\mathbb A_{\mathbb X_1}^*\otimes \cdots\otimes \mathbb A_{\mathbb X_n}^*)^*$$
\end{proposition}

\begin{proof} It is a consequence of the equalities
$$\aligned & \mathbb Hom(\mathbb X_1\times \cdots\times \mathbb X_n,\mathbb M^*)=
\mathbb Hom(\mathbb X_1,\mathbb Hom(\mathbb X_2\times \cdots\times \mathbb X_n,\mathbb M^*))\\ & \overset{\text{Induction}}=\mathbb Hom(\mathbb X,\mathbb Hom_{\mathcal R}(\mathbb B_{\mathbb X_2}\otimes \cdots\otimes \mathbb B_{\mathbb X_n},\mathbb M^*))\\ & \overset{\text{\ref{2.3}}}=\mathbb Hom_{\mathcal R}(\mathbb B_{\mathbb X_1},\mathbb Hom_{\mathcal R}(\mathbb B_{\mathbb X_2}\otimes \cdots\otimes \mathbb B_{\mathbb X_n},\mathbb M^*))=\mathbb Hom_{\mathcal R}(\mathbb B_{\mathbb X_1}\otimes \cdots\otimes \mathbb B_{\mathbb X_n},\mathbb M^*)\endaligned$$
\end{proof}

\begin{proposition} Let $\mathbb X,\mathbb Y$ be affine functors such that
$\mathbb A_{\mathbb X}, \mathbb A_{\mathbb Y}\in\mathfrak F$, then
$\mathbb X\times\mathbb Y$ is an affine functor and  $\mathbb
A_{\mathbb X\times \mathbb Y}\in\mathfrak F$.\end{proposition}

\begin{proof} $\mathbb A_{\mathbb X\times \mathbb Y}=(\mathbb A_{\mathbb X}^*\otimes \mathbb A_{\mathbb Y}^*)^*=\mathbb A_{\mathbb X}\tilde\otimes \mathbb A_{\mathbb Y}
\in\mathfrak F$  and
$\Spec \mathbb A_{\mathbb X\times \mathbb Y}=
\Spec (\mathbb A_{\mathbb X}\otimes \mathbb A_{\mathbb Y})=\mathbb X\times \mathbb Y$,  by Theorem \ref{n5.22}.

\end{proof}

\begin{proposition} \label{prodafin} If $\mathbb X=\mathbb X_1\times \mathbb X_2$ is an affine functor and $\mathbb X(R)\neq \emptyset$ then $\mathbb X_1$ is an affine functor.\end{proposition}

\begin{proof} Given $(x_1,x_2)\in \mathbb X_1(R)\times \mathbb X_2(R)=\mathbb X(R)$, let
$i\colon \mathbb X_1\hookrightarrow \mathbb X$, $i(x)=(x,x_2)$. Let $\pi_1\colon
\mathbb X =\mathbb X_1\times \mathbb X_2\to \mathbb X_1$, $\pi_1((y_1,y_2))=y_1$.
Let $i^*\colon \mathbb A_{\mathbb X} \to \mathbb A_{\mathbb X_1}$, $\pi_1^*\colon
\mathbb A_{\mathbb X_1}\to \mathbb A_{\mathbb X}$ be the morphisms induced by $i$ and $\pi_1$ respectively. Obviously, $\pi_1^*\circ i^*=Id$ because $i\circ \pi_1=Id$.
Hence, $\mathbb A_{\mathbb X_1}$ is  a direct summand of $\mathbb A_{\mathbb X}$, and
it is a reflexive functor because  $\mathbb A_{\mathbb X}$ is reflexive.

Given, $f\in \mathbb Hom_{\mathcal R-alg}(\mathbb A_{\mathbb X_1},\mathcal R)$
then $f\circ i^*=(y_1,y_2)\in \mathbb Hom_{\mathcal R-alg}(\mathbb A_{\mathbb X},\mathcal R)=\mathbb X$ and $f=f\circ i^*\circ \pi_1^*=(y_1,y_2)\circ\pi_1=y_1$, that is,
the morphism $\mathbb X_1\to {\rm Spec}\, \mathbb A_{\mathbb X_1}$ is surjective. Finally, since the composition
$$\mathbb X_1\to {\rm Spec}\, \mathbb A_{\mathbb X_1}\overset{i^{**}}\to
{\rm Spec}\, \mathbb A_{\mathbb X}=\mathbb X$$
is equal to the morphism $i$, then $\mathbb X_1={\rm Spec}\, \mathbb A_{\mathbb X_1}$.

\end{proof}

\section{Formal schemes}

\begin{definition} Let $\mathcal C^*\in\mathfrak F$ be a scheme of commutative algebras.
We will say that
$\Spec\mathcal C^*$ is a formal scheme.
If $\Spec\mathcal C^*$ is a functor of monoids we will say that it is a formal monoid. \end{definition}

Recall that if $C$ is a free module $\mathcal C^*\in\mathfrak F$.

\begin{note} \label{n3.22} By Note \ref{n3.2}, formal schemes
are affine functors. In fact,   $\Spec\mathcal C^*$ is a direct limit of finite $R$-schemes. Reciprocally, if $R$ is a field, a direct limit of finite $R$-schemes are formal schemes, by Theorem \ref{4.4}.
If $R$ is a field, Demazure (\cite{Demazure2}) defines a formal scheme as a functor (over the $R$-finite dimensional rings) which is a direct limit of finite $R$-schemes.
\end{note}

The direct product $\Spec\mathcal C_1^*\times \Spec \mathcal C_2^*=\Spec (\mathcal C_1^*\tilde \otimes\mathcal C_2^*)=\Spec (\mathcal C_1 \otimes\mathcal C_2)^* $ of formal schemes is a formal scheme.

\begin{example}
Let $X$ be a set. Let us consider the discrete topology on $X$.
Let $\bf X$ be ``the constant
functor $X$'', defined by $${\bf X}(S) := {\Aplic}
({\rm Spec}\, S, X)$$ for every commutative $R$-algebra $S$. If ${\rm Spec}\,S$ is connected then ${\bf X}(S)=X$.

Let ${\mathbb A}_X$ be the functor of algebras defined by
$${\mathbb A}_X(S) :=  \prod_X S$$
for each commutative $R$-algebra $S$. Observe that
${\mathbb A}_X = \underset{X}{\prod} {\mathcal R}
=(\underset{X}{\oplus} {\mathcal R})^*$ is a commutative algebra
scheme. $\bf X$ is a formal scheme because ${\rm Spec}\, {\mathbb A}_X = {\bf X}$:

$$\aligned ({\rm Spec}\, {\mathbb
A}_X )(S) & = {\rm Hom}_{\mathcal R-alg} (\prod_X\mathcal R, {\mathcal S} )
\overset{\text{\ref{a1.8}}}=\ilim{\underset{|Y|<\infty}{Y\subset X}}{\rm Hom}_{\mathcal R-alg} (\prod_Y\mathcal R, {\mathcal S} )\\ &=\ilim{\underset{|Y|<\infty}{Y\subset X}}
{\rm Hom}_{R-alg} (\prod_Y R, {S} )=\ilim{\underset{|Y|<\infty}{Y\subset X}}{\Aplic}(\Spec S,Y)\\ & ={\Aplic}(\Spec S,X)={\bf X}(S)\endaligned$$

\end{example}

Obviously, ${\rm Spec}\,  (\underset{\underset{i\in I}{\longrightarrow}}{\lim}\,
{\mathbb A}_i) =  \underset{\underset{i\in I}{\longleftarrow}}{\lim}\, ({\rm Spec }\, {\mathbb A_i})$.

\begin{theorem} \label{4.4} Let $\{{\Spec \mathcal C_i^*}\}_{i\in I}$ be a direct system of formal schemes. Then, $$\ilim{i} \Spec \mathcal C_i^* =
\Spec (\ilim{i} \mathcal C_i)^*$$ and it is an affine functor.\end{theorem}

\begin{proof} Write $C=\ilim{i\in I} C_i$, then $\mathcal C^*=\plim{i\in I} \mathcal C_i^*$.
$ {\rm Hom}_{\mathcal R}( {\mathcal C^*}, {\mathcal S})
={\rm Hom}_{\mathcal R}( (\ilim{i\in I} {\mathcal C_i})^*, {\mathcal S})
=(\ilim{i\in I} {C_i})\otimes S=
\underset{\underset{i\in I}{\longrightarrow}}{\lim}\,
{\rm Hom}_{\mathcal R}( {\mathcal C_i}^*,{\mathcal S})$.
Likely, ${\rm Hom}_{\mathcal R}( {\mathcal C^*\otimes \mathcal C^*}, {\mathcal S})=
{\rm Hom}_{\mathcal R}( (\mathcal C\otimes \mathcal C)^*, {\mathcal S})=
\underset{\underset{i\in I}{\longrightarrow}}{\lim}\,
{\rm Hom}_{\mathcal R}( ({\mathcal C_i}\otimes \mathcal C_i)^*,{\mathcal S})=
\underset{\underset{i\in I}{\longrightarrow}}{\lim}\,
{\rm Hom}_{\mathcal R}( {\mathcal C_i}^*\otimes {\mathcal C_i}^*,{\mathcal S})$.
Then the kernel of the morphism $ {\rm Hom}_{\mathcal R}( {\mathcal C^*}, {\mathcal S})\to {\rm Hom}_{\mathcal R}( {\mathcal C^*\otimes \mathcal C^*}, {\mathcal S})$, $f\mapsto \tilde f$, where $\tilde f(c_1\otimes c_2)=f(c_1c_2)-f(c_1)f(c_2)$ coincides with the kernel of the morphism
$\underset{\underset{i\in I}{\longrightarrow}}{\lim}\,
{\rm Hom}_{\mathcal R}( {\mathcal C_i}^*,{\mathcal S})\to
\underset{\underset{i\in I}{\longrightarrow}}{\lim}\,
{\rm Hom}_{\mathcal R}( {\mathcal C_i}^*\otimes {\mathcal C_i}^*,{\mathcal S})$, $(f_i)\mapsto (\tilde f_i)$.
Then, ${\rm Hom}_{\mathcal R-alg}(  {\mathcal C^*}, {\mathcal S})  =
\underset{\underset{i\in I}{\longrightarrow}}{\lim}\,
{\rm Hom}_{\mathcal R-alg}( {\mathcal C_i}^*,{\mathcal S})$ and

$$(\Spec {\mathcal C^*})
 (S)=(\underset{\underset{i\in I}{\longrightarrow}}{\lim}\,\Spec {\mathcal C_i}^*)(S)$$
Finally, $\mathbb A_{\Spec {\mathcal C^*}}=\plim{i} \mathbb A_{\Spec {\mathcal C_i}^*}=
\plim{i} \mathcal C_i^*=\mathcal C^*$.

\end{proof}

From now on, in this section, we will assume that $R=K$ is a field. \label{seccion4}

\begin{definition} \label{4.1} Let $X$ be a $K$-scheme and $I$ the set of finite subschemes of $X$.
Given $K$-scheme $Y$ write $A_Y :=\mathcal O_Y(Y)$, the ring of functions of $Y$.
Define $\bar {\mathcal A}_X:=\plim{i\in I} \mathcal A_i$ and $$\bar X:=\Spec \bar{\mathcal A}_X\overset{\text{\ref{4.4}}}=\ilim{i\in I} \Spec \mathcal A_i$$ That is, ``$\bar X$ is the direct limit of the set of all finite subschemes of $X$''.
\end{definition}

We have a natural morphism $\bar X\hookrightarrow X$.

\begin{theorem} \label{universal} Let $X$ be a $K$-scheme. Then: $$\Hom(\Spec {\mathcal C^*}, X)=
\Hom(\Spec {\mathcal C^*}, \bar X)$$ for every formal scheme $\Spec\mathcal C^*$.\end{theorem}

\begin{proof} ${\mathcal C}^*=\plim{i} \mathcal S_i$, where the $S_i$ are finite $K$-algebras.
Then,
 $$\begin{array}{l} \Hom(\Spec {\mathcal C^*}, X)=\Hom(\ilim{i} \Spec {\mathcal S_i}, X) =
\plim{i}\Hom( \Spec {\mathcal S_i}, X)\\= \plim{i} \Hom(\Spec {\mathcal S_i}, \bar X)= \Hom(\ilim{i} \Spec {\mathcal S_i}, \bar X)=
\Hom(\Spec {\mathcal C^*}, \bar X)\end{array} $$\end{proof}

Let $X$ and $Y$ be two $K$-schemes. By the universal property \ref{universal}, it can be checked that
$$\overline{X\times Y}=\bar X\times \bar Y$$

For every functor of $\mathcal K$-modules ${\mathbb M}$ there exists a natural morphism from the quasi-coherent module associated with  ${\mathbb M}(K)$ into ${\mathbb M}$. If $\mathcal
C^*$ is a commutative algebra scheme, then we have a natural morphism from the quasi-coherent algebra associated with $C^*$ into $\mathcal C^*$ and therefore a natural morphism $\Spec{\mathcal C}^*\to \Spec C^*$. Moreover, this is injective, because $\mathcal C^*=\plim{i} \mathcal C_i$, (where $\mathcal C_i$ are the quotients of coherent algebras of $\mathcal C^*$) and:

$$\aligned (\Spec \mathcal C^*)(S) & =\Hom_{\mathcal K-alg}({\mathcal C^*}, {\mathcal S})=\ilim{i} \Hom_{\mathcal K-alg}({\mathcal C_i}, {\mathcal S})\\ &=\ilim{i} \Hom_{K-alg}({C_i}, {S})\subseteq \Hom_{K-alg}({C^*}, {S})=(\Spec  C^*)(S)\endaligned$$

Let $X$ be a compact and separated $K$-scheme and $\bar A_X=\bar{\mathcal A}_X(K)=\plim{i} A_i$. We can define a natural morphism
$\Spec  \bar A_X \to X$:

 Let $X':=\{x\in X\colon $ the residual field of $x$ is a finite extension of $K$, $\dim_K \mathcal O_{X,x}/{\mathfrak p}_x<\infty\}$. For any subset $J\subset X'$ let us denote
$\bar A_J:=\prod_{x\in J} \hat{\mathcal O}_{X,x}$,
where $$\hat{\mathcal O}_{X,x}:=\plim{ \dim_K\mathcal O_{X,x}/I<\infty} \mathcal O_{X,x}/I$$
Write $\bar J :=\Spec \bar A_J$. For any two disjoint subsets $J,J'\subseteq X'$, $\bar A_{J\coprod J'}=\bar A_J\times \bar A_{J'}$, and, in general, $\bar A_{J\cup J'}=\bar A_{J\cap J'}\times
\bar A_{J-J\cap J'}\times \bar A_{J'-J\cap J'}$. $\bar J$ is an open and closed subset of $\overline{X'}$
and we have that $\overline{J\cup J'}=\bar J\cup \overline{J'}$ and $\overline{J\cap J'}=\bar J\cap \overline{J'}$.

For any affine open subset $U=\Spec A\subset X$ we have a natural morphism $A\to \prod_{x\in U'} \hat{\mathcal O}_{X,x}=\bar A_U$, and therefore a natural morphism
$f_U\colon  \overline{U'}\to U\subset X$. For any other affine open subset $V\subset X$ we have another morphism $f_V\colon \overline{V'}\to V\subset X$. As $U\cap V$ is affine and $f_U$ is equal to $f_V$ on $\overline{(U\cap V)'}=\overline{U'}\cap \overline{V'}$,  we have a well-defined natural morphism $\Spec  \bar A_X=\overline{X'} \to X$.

The natural morphism $\bar X\to X$ is equal to the composition of the natural morphisms:
$\bar X=\Spec   \bar{\mathcal  A}_X\to \Spec  \bar A_X \to X$, because for any finite subscheme $X_i$ of $X$ the composition morphism $X_i\to \Spec  \bar A_X \to X$ is equal to the natural inclusion $X_i\hookrightarrow X$.

\begin{theorem} \label{8} Every morphism $\Spec \mathcal C^*\to X$, from a formal scheme to a compact and separated $K$-scheme uniquely factors via $\Spec C^*$, that is,

$$\Hom(\Spec {\mathcal C}^*, X)=\Hom_{K-sch}(\Spec {C}^*, X).$$ \end{theorem}

\begin{proof}  If $X=\Spec A$ is an affine scheme then
$$\aligned \Hom(\Spec\mathcal C^*,\Spec A) & =\Hom_{\mathcal R-alg}(\mathcal A, \mathcal C^*)=
\Hom_{\mathcal R-alg}(A, C^*)\\ & =\Hom_{K-sch}(\Spec C^*,\Spec A)\endaligned$$

By Theorem \ref{universal}, every morphism from $\Spec {\mathcal C}^*$ into $X$ uniquely factors through
$\bar X$. The morphism $\bar X\to X$ factors through $\Spec \bar A_{X}$ and the morphism  $\Spec {\mathcal C}^*\to \Spec  \bar A_{X}$
uniquely factors through $\Spec {C}^*$. So that we have the commutative diagram:
$$\xymatrix{ \Spec {\mathcal C}^* \ar[r] \ar[d] & \Spec C^* \ar@{-->}[d] & \\ \bar X \ar[r] & \Spec  \bar A_{X} \ar[r] & X}$$
Then, the natural morphism $\Hom(\Spec {C}^*, X)\to \Hom(\Spec {\mathcal C}^*, X)$ is surjective.

Let us write ${\mathcal C^*}=\plim{i}\mathcal C_i$, with $\dim_K C_i<\infty$ and $\mathcal C^*\to \mathcal C_i$ surjective. The smallest closed subscheme  of $\Spec C^*$ containing every $\Spec C_i$ is
$\Spec C^*$. Therefore,

$$\aligned \Hom_{K-sch}(\Spec {C}^*, X) & \subseteq \plim{i}
\Hom_{K-sch}( \Spec {C_i}, X) =\Hom(\ilim{i} \Spec {\mathcal C_i}, X)\\ & = \Hom(\Spec {\mathcal C}^*, X)\endaligned$$

\end{proof}

\begin{note} \label{3.9} Let $X$ and $Y$ be compact and separated $K$-schemes. Every commutative diagram
$$\xymatrix{\Spec{\mathcal B^*} \ar[r] \ar[d] & X \ar[d] \\
\Spec{\mathcal C^*} \ar[r] & Y}$$ uniquely extends to a commutative diagram

$$\xymatrix{\Spec{\mathcal B^*} \ar[r] \ar[d] & \Spec B^* \ar[r] \ar@{-->}[d] &  X \ar[d]\\
\Spec{\mathcal C^*} \ar[r] & \Spec C^* \ar[r] & Y}$$ because the composition  $\Spec{\mathcal B^*}\to \Spec{\mathcal C^*}\to \Spec C^*$ factors through $\Spec B^*$.

\end{note}

\section{Completion along a closed subscheme}

\begin{proposition} \label{prodq} Let $M_{n+1}\to M_n$ be epimorphisms of $R$-modules, for all $n\in \mathbb N$ and let  $N$ be an  $R$-module. It holds that
$$\Hom_{\mathcal R}(\plim{n\in\mathbb N} \mathcal M_n,\mathcal N)=
\ilim{n\in\mathbb N} \Hom_{\mathcal R}(\mathcal M_n,\mathcal N)$$
Hence $(\plim{n\in\mathbb N} \mathcal M_n)^*=\ilim{n\in\mathbb N} \mathcal M_n^*$ and $\plim{n\in\mathbb N} \mathcal M_n$ is a reflexive $\mathcal R$-module.\end{proposition}

\begin{proof} Let $f\in \Hom_{\mathcal R}(\plim{n\in\mathbb N} \mathcal M_n,\mathcal N)$. Firstly, let us prove that the morphism
$f_R\colon \plim{n\in\mathbb N} M_n\to  N$ induced by $f$ factors through $M_r$, for some
$r\in\mathbb N$: suppose that, for each $r$, there exists an element $s_r=(m_n)\in \plim{n\in\mathbb N}  M_n\subset \prod_n M_n$, such that $m_r= 0$ and $f_R(s_r)\neq 0$.
The morphism $g\colon \prod_{n\in\mathbb N}\mathcal R\to \mathcal N$,
$g((a_n)_n):=f(\sum_n a_n\cdot s_n)$ satisfies that $g_{|\mathcal R}\neq 0$ for every factor $\mathcal R\subset \prod_n\mathcal R$ and this contradicts the fact that ${\rm Hom}_{\mathcal R}(\prod_n \mathcal R, \mathcal N)\overset{\text{\ref{a1.8}}}=(\oplus_n R)\otimes N=\oplus_n{\rm Hom}_{\mathcal R}(\mathcal R,\mathcal N)$. Then, $f_R$ factors through a (unique) morphism $h_R\colon M_r\to N$, for some $r$.

Next, given a commutative $R$-algebra $S$, let us check that the  morphism
$$f_S\colon \plim{n\in\mathbb N}  (M_n\otimes_RS)\to  N\otimes_RS$$ induced by $f$ factors through $h_S\colon M_r\otimes_RS\to N\otimes S$, $h_S(m_r\otimes s):=h_R(m_r)\otimes s$: there exists $r'\geq r$ such that
$f_S$ factors through a morphism $h'\colon M_{r'}\otimes_R S\to N\otimes S$. Given $m_{r'}\in M_{r'}$, let $(m'_n)\in \plim{n\in\mathbb N} \mathcal M_n$ such that
$m_{r'}=m'_{r'}$. Then, $h'(m_{r'}\otimes 1)=f_S((m'_{n}\otimes 1))=f_R((m'_{n}))\otimes 1=h_R(m'_{r})\otimes 1=h_S(m'_{r}\otimes 1)$ and $h'$ factors through $h_S$.
Hence $f_S$ factors through $h_S$.

\end{proof}

\begin{definition} \label{eje36} Let $A$ be a commutative $R$-algebra and $I\subset A$ an ideal. Let $\hat{\mathcal A}:=\plim{n\in\mathbb N} \mathcal A/\mathcal I^n$. We will say that $\Spec\hat{\mathcal A}$ is the completion of
$\Spec A$ along the closed set $\Spec A/I$.\end{definition}

$\hat A$ is a reflexive functors of $R$-modules by Proposition \ref{prodq}. $\hat{\mathcal A}$ is proquasi-coherent algebra: $\hat{\mathcal A}\,^*=\ilim{n} (\mathcal A/\mathcal I^n)^*$, then every morphism of functors of algebras $\hat{\mathcal A}\to \mathcal B$ factors through some  $\mathcal A/\mathcal I^n$. Hence the inverse limit of the proquasi-coherent algebra quotients of $\hat{\mathcal A}$ is equal to $\plim{n\in\mathbb N} \mathcal A/\mathcal I^n=
\hat{\mathcal A}$.

\begin{proposition} \label{eje4.6} The completion of
$\Spec A$ along the closed set $\Spec A/I$ is an affine functor.\end{proposition}

\begin{proof}
$\Hom_{\mathcal R-alg}(\hat{\mathcal A}, \mathcal C)=
\ilim{n\in\mathbb N}\Hom_{\mathcal R-alg}(\mathcal A/\mathcal I^n, \mathcal C)$, then $\Spec \hat{\mathcal A} = \ilim{n\in\mathbb N}
\Spec  \mathcal A/\mathcal I^n.$
By \ref{3.10}, $\mathbb A_{\Spec \hat{\mathcal A}}=\plim{n} \mathbb A_{\Spec \mathcal A/\mathcal I^n}=\plim{n} \mathcal A/\mathcal I^n= \hat{\mathcal A}$.
\end{proof}

Let $B$ be a commutative $R$-algebra, $J\subset B$ be an ideal and $\hat{\mathcal B}=\plim{n\in\mathbb N} \mathcal B/\mathcal J^n$.
Consider the ideal $I\otimes B+A\otimes J\subseteq A\otimes B$. Then,

$$\aligned  \Spec & \widehat{\mathcal A \otimes\mathcal B} \! =\!\mathbb Hom_{\mathcal R-alg}(\widehat{\mathcal A\otimes\mathcal B},\mathcal K)\!=\!\ilim{n} \mathbb Hom_{\mathcal R-alg}(\mathcal A\otimes\mathcal B/(\mathcal I\otimes\mathcal B+\mathcal A\otimes \mathcal J)^n,\mathcal K)
\\ &
= \ilim{m} \mathbb Hom_{\mathcal R-alg}(\mathcal A/\mathcal I^m\otimes \mathcal B/\mathcal J^m,\mathcal K)\\ & =\ilim{m} ( \mathbb Hom_{\mathcal R-alg}(\mathcal A/\mathcal I^m,\mathcal K)\times
\mathbb Hom_{\mathcal R-alg}(\mathcal B/\mathcal J^m,\mathcal K))=
\Spec \hat{\mathcal A}\times \Spec \hat{\mathcal B}\endaligned$$
and
$$\hat{\mathcal A}\tilde \otimes \hat{\mathcal B}=(\hat{\mathcal A}^*\otimes \hat{\mathcal B}^*)^*=(\ilim{n} ((\mathcal A/\mathcal I^n)^*\otimes (\mathcal B/\mathcal J^n)^*))^*
\overset{\text{\ref{2.21}}}=
\plim{n} (\mathcal A/\mathcal I^n\otimes \mathcal B/\mathcal J^n)
=\widehat{\mathcal A\otimes\mathcal B}$$




\section{Affine functors of monoids}

Let $\mathbb G$ be a functor  of monoids. ${\mathcal R}\mathbb G$ is obviously a functor of
${\mathcal R}$-algebras. Given a functor of ${\mathcal R}$-algebras $\mathbb B$, it is
easy to check the equality $${\rm Hom}_{mon} (\mathbb G, {\mathbb B})
= {\rm Hom}_{{\mathcal R}-alg} ({\mathcal R}\mathbb G, \mathbb B) .$$
The closure of dual functors of $\mathcal R$-algebras of
$\mathbb G$ is equal to the closure of dual functors of $\mathcal R$-algebras of $\mathcal R\mathbb G$.

\begin{theorem}\label{2.5}
Let $\mathbb G$ be a functor of monoids with a reflexive functor of functions. Then, the
closure of dual functors of algebras of $\mathbb G$ is
$\mathbb A_{\mathbb G}^*$. That is,

$${\rm Hom}_{mon}(\mathbb G, {\mathbb B}) = {\rm Hom}_{\mathcal R-alg}(\mathcal R\mathbb G,
{\mathbb B}) = {\rm Hom}_{\mathcal R-alg}({\mathbb A_{\mathbb G}^*}, {\mathbb B})$$
for every dual functor of $\mathcal R$-algebras $\mathbb B$.
\end{theorem}

\begin{proof}
$(\mathcal R\mathbb G)^*=\mathbb A_{\mathbb G}$ is reflexive, so the closure of dual functors of
algebras of $\mathbb G$ is $\mathbb A_{\mathbb G}^*$, by Proposition \ref{n2.16}.
\end{proof}

\begin{theorem}\label{Gmodulos}
Let $\mathbb G$ be a functor of monoids with a reflexive functor of functions.  The category of quasi-coherent ${\mathbb G}$-modules is equi\-valent to the
category of quasi-coherent $\mathbb A_{\mathbb G}^*$-modules.

Likewise, the category of
dual functors of $\mathbb G$-modules is equivalent to the category of dual
functors of $\mathbb A_{\mathbb G}^*$-modules.
\end{theorem}

\begin{proof} It a consequence of Proposition \ref{n2.16}.
\end{proof}

Remark that the structure of functor of algebras of
$\mathbb A_{\mathbb G}^*$ is the only one that makes the morphism ${\mathbb G}\to \mathbb A_{\mathbb G}^*$ a morphism of functors of
monoids.

\begin{theorem}
Let $\mathbb G$ be a functor of monoids with a reflexive functor of functions and let
$\mathbb M,\mathbb M'$ be reflexive functors of $\mathbb G$-modules.
Assume that $\mathbb A_{\mathbb G}, \mathbb M,\mathbb M'\in\mathfrak F$ or that $R=K$ is a field. Then, a morphism
    of $\mathcal R$-modules $\mathbb M\to\mathbb M'$ is a morphism of $\mathbb G$-modules if and only if $\mathbb M(R)\to\mathbb M'(R)$ is a morphism of $\mathbb A^*_{\mathbb G}(R)$-modules. Let  $\mathcal M$ be a $\mathbb G$-module, then the set of quasi-coherent $\mathbb G$-submodules of $\mathcal M$ is equal to the set of $\mathbb A_{\mathbb G}^*(R)$-submodules of $M$.

\end{theorem}

\begin{proof} It a consequence of Theorem \ref{Gmodulos} and Proposition \ref{nS1}.
\end{proof}

\begin{example}
The ${\mathbb C}$-linear representations of $(\mathbb Z,+)$ are equivalent to
the $\mathbb A_{\mathbb Z}^*$-modu\-les. $\mathbb A_{\mathbb Z}^*$ is the quasi-coherent algebra associated with
${\mathbb
C}[x,1/x]$. The category of $\mathbb A_{\mathbb Z}^*$-modu\-les is equal to the category of ${\mathbb
C}[x,1/x]$-modules.
The natural morphism $\mathbb Z\to \mathbb A_{\mathbb Z}^*$ assigns $n$ to $x^n$. Thus, if
$V$ is a finite ${\mathbb C}$-linear representation of $\mathbb Z$, then
$$V= \underset{\alpha,n,m}{\oplus} ({\mathbb
C}[x]/(x-\alpha)^{n})^{m},\qquad (\alpha\neq 0)$$ such that $r\cdot (\overline{
p_{\alpha,n,m}(x)})_{\alpha,n,m}=(\overline{x^r\cdot
p_{\alpha,n,m}(x)})_{\alpha,n,m}$.
\end{example}

Let $\mathbb G$ be a functor of monoids with a reflexive functor of functions. Let $m\colon \mathbb G\times \mathbb G\to \mathbb G$ be the multiplication morphism. Then, the composition morphism
of $m$ with the natural morphism $\mathbb G\to \mathbb A_{\mathbb G}^*$ factors through $\mathbb A^*_{\mathbb G}\otimes \mathbb A^*_{\mathbb G}$, by \ref{2.3} and \ref{244}, that is, we have a commutative diagram
$$\xymatrix{\mathbb G\times \mathbb G \ar[r]^-m \ar[d]& \mathbb G \ar[d]\\ \mathbb A^*_{\mathbb G}\otimes \mathbb A^*_{\mathbb G} \ar[r]^-m & \mathbb A^*_{\mathbb G}}$$
Let $e\in \mathbb G(R)\subset \mathbb A^*_{\mathbb G}(R)$ the unit  of $\mathbb G$. Then, we can define a morphism $e\colon \mathcal R\to \mathbb A^*_{\mathbb G}$. It is easy to check that $\{\mathbb A^*_{\mathbb G},m,e\}$ is a functor of $\mathcal R$-algebras. Moreover, the dual morphisms of the multiplication morphism $m$ and the unit morphism $e$ are
the natural morphisms $\mathbb A_{\mathbb G}\to \mathbb A_{\mathbb G\times \mathbb G}$ and
$\mathbb A_{\mathbb G}\overset e\to \mathcal R$, which
are morphisms of $\mathcal R$-algebras.

Conversely, let $\mathbb X$ be an affine functor and
assume that $\mathbb A^*_{\mathbb X}$ is a functor of $\mathcal R$-algebras, such that
the dual morphisms $m^*$ and $e^*$, of the multiplication morphism $m\colon \mathbb A_{\mathbb X}^*\otimes
\mathbb A_{\mathbb X}^*\to \mathbb A_{\mathbb X}^*$ and the unit morphism $e\colon \mathcal R\to \mathbb A_{\mathbb X}^*$ are morphisms of $\mathcal R$-algebras. Given a point
$(x,x')\in \mathbb X\times \mathbb X\subset \mathbb Hom_{\mathcal R-alg}(\mathbb A_{\mathbb X\times\mathbb X},\mathcal R)$ then $(x,x')\circ m^*\in \mathbb Hom_{\mathcal R-alg}(\mathbb A_{\mathbb X},\mathcal R)=\mathbb X$ and we have the commutative diagram
$$\xymatrix{\mathbb X\times \mathbb X \ar@{-->}[r]^-m \ar[d]& \mathbb X \ar[d]\\ \mathbb A^*_{\mathbb X}\otimes \mathbb A^*_{\mathbb X} \ar[r]^-m & \mathbb A^*_{\mathbb X}}$$
Obviously $e\in \mathbb Hom_{\mathcal R-alg}(\mathbb A_{\mathbb X},\mathcal R)=\mathbb X$.
It is easy to check that $\{\mathbb X,m,e\}$ is a functor of monoids.

\begin{definition} An affine functor $\mathbb G=\Spec \mathbb A$ is said to be an
affine functor of monoids  if $\mathbb G$ is a functor of monoids.
\end{definition}

\begin{example} Affine $R$-monoid schemes, formal monoids,  the completion of an affine monoid scheme along a closed submonoid scheme, $\mathcal V$, $\mathbb End_{\mathcal R}\mathcal V$ ($V$ being a free $R$-module) are examples of affine functors of monoids, by \ref{eje4.5}, \ref{eje4.6} and \ref{pro314}.
\end{example}

Let $\mathbb G$ and $\mathbb G'$ be affine functors of monoids. Then,
$$\Hom_{mon}(\mathbb G,\mathbb G')=\{f\in \Hom_{\mathcal R}(\mathbb A_{\mathbb G'},\mathbb A_{\mathbb G})\colon f,\,f^*\text{ are morph. of funct. of $\mathcal R$-alg.}\}:$$
Let $h\colon \mathbb G\to \mathbb G'$ be a morphism of functors of monoids. The composition
morphism of $h$ with the natural morphism $\mathbb G'\to \mathbb A_{\mathbb G'}^*$ factors through $\mathbb A_{\mathbb G}^*$, that is, we have a commutative diagram
$$\xymatrix{\mathbb G\ar[r]^-h \ar[d]& \mathbb G' \ar[d]\\ \mathbb A^*_{\mathbb G}\ar[r] & \mathbb A^*_{\mathbb G'}}$$
The dual morphism $\mathbb A_{\mathbb G'}\to \mathbb A_{\mathbb G}$ is the morphism induced by $h$ between the functors of functions. Inversely, let $f\colon \mathbb A_{\mathbb G'}\to \mathbb A_{\mathbb G}$ be a morphism of functors of $\mathcal R$-algebras, such that $f^*$ is also a morphism  of functors of $\mathcal R$-algebras. Given $g\in\mathbb G$, then $f^*(g)=g\circ f\in
\mathbb Hom_{\mathcal R-alg}(\mathbb A_{\mathbb G'},\mathcal R)=\mathbb G'$.
Hence, $f^*_{|\mathbb G}\colon \mathbb G\to \mathbb G'$ is a morphism of functors of monoids.

\begin{theorem} \label{5.3n} Let $\mathbb G=\Spec \mathbb A_{\mathbb G}$ be an affine functor.
Giving to $\mathbb G$ a structure of functor of monoids is equivalent to giving
to $ \mathbb A_{\mathbb G}$ a structure of functor of bialgebras. Let $\mathbb G$ and $\mathbb G'$ be two affine functors of monoids, it holds that
$$\Hom_{mon}(\mathbb G,\mathbb G')=\Hom_{bialg}( \mathbb A_{\mathbb G'},\mathbb A_{\mathbb G})$$

\end{theorem}

\begin{theorem} \label{5.3m} The category of cocommutative bialgebras $A$ is equivalent to the category of formal monoids $\Spec \mathcal A^*$ (we assume the $R$-modules $A$ are projective).\end{theorem}


\begin{theorem} \label{4.4b} Let $G$ be a $K$-scheme on groups (resp. monoids). Then $\bar G$ is a functor of groups (resp. monoids), the natural morphism $\bar G\to
G$ is a morphism of functors of monoids and
$$\Hom_{mon}(\Spec {\mathcal C^*}, G)=
\Hom_{mon}(\Spec {\mathcal C^*}, \bar G)$$ for every formal  monoid $\Spec\mathcal C^*$. If $G$ is commutative, then $\bar G$ is commutative.\end{theorem}

\begin{proof} Let $\mu\colon G\times G\to G$ the multiplication morphism.
By Theorem \ref{universal}, the composition morphism $\overline{G\times G}=\bar G\times \bar G\to G\times G\to G$ factors through a unique morphism $\mu'\colon \bar G\times \bar G\to \bar G$, that is, we have the commutative diagram:
$$\xymatrix{\bar G\times \bar G \ar[r] \ar[d]^-{\mu'} &
G \times G \ar[d]^-\mu \\ \bar G \ar[r] & G}$$
Let $*\colon G\to G$ be the inverse morphism. The composition
$\bar G \to G \overset{*}\to G$ factors through a unique morphism
$*'\colon \bar G\to  \bar G$, that is, we have the commutative diagram:
$$\xymatrix{ \bar G \ar[r] \ar[d]^-{*'} & G \ar[d]^-{*} \\
\bar G\ar[r] & G}$$
Now it is easy to check that $(\bar G,\mu',*')$ is a functor of groups and to conclude the proof.\end{proof}

\begin{proposition} \label{1.3} Let $\Spec {\mathcal C^*}$ be a formal monoid (resp. group) and let  $G$ be a compact $K$-scheme on monoids. Let  $f\colon \Spec{\mathcal
C^*}\to G$ be a morphism of functors of monoids and let $f'\colon \Spec
C^*\to G$ be the induced morphism.

Then the (closed) scheme-theoretic image of $f'$, ${\rm Im}\, f'$,  is a $K$-subscheme on monoids (resp. groups) of
$G$. If $\Spec {\mathcal C^*}$ is an abelian formal monoid, then
${\rm Im}\, f'$ is abelian.\end{proposition}

\begin{proof} Let $\mu'\colon \Spec {\mathcal C^*}\times \Spec {\mathcal C^*}\to \Spec
{\mathcal C^*}$ and $\mu\colon G\times G\to
G$ be the operations of the monoids. Consider the commutative diagram:
$$\xymatrix{\Spec{\mathcal C^*}\times \Spec{\mathcal C^*} \ar[r]
\ar[d]^-{\mu'} &  \Spec (C\otimes C)^* \ar[d] \ar[r] &
\Spec{C^*}\times \Spec{C^*} \ar[r]^-{f'\times f'}  & G\times G
\ar[d]^-\mu \\ \Spec{\mathcal C^*} \ar[r] & \Spec{ C^*} \ar[rr]^-{f'}&
& G}$$

The scheme-theoretic image of $\Spec (C\otimes C)^*$ in $G\times
G$ is equal to ${\rm Im}\,(f'\times f')$. As ${\rm Im}\,(f'\times f')={\rm Im}\,
f'\times {\rm Im}\, f'$, we have the commutative diagram:
$$\xymatrix{{\rm Im}\, f'\times {\rm Im}\, f' \ar@{^{(}->}[r] \ar[d] & G\times G \ar[d] \\ {\rm Im}\,
f' \ar@{^{(}->}[r] & G} $$

If we denote with $*$ the inverse morphism  (with respect to the group law), then we have the commutative square:
$$\xymatrix{\Spec {\mathcal C^*} \ar[d]^-{*} \ar[r] & \Spec C^* \ar[r]^-{f'} \ar[d] & G \ar[d]^-{*} \\ \Spec {\mathcal C^*} \ar[r] & \Spec C^* \ar[r]^-{f'} & G }$$ So the inverse of $G$ restricts to ${\rm Im}\, f'$. The rest of the details can be checked in a similar way.

\end{proof}

\begin{note}
 Let $X$ be a $K$-scheme and $x\colon {\rm Spec}\, K\to X$ a rational point. Denote by $\hat X$ the direct limit of finite subschemes  $X_i\subset X$ with support on $x$. Let $\hat {\mathcal A}_X:=\underset{\leftarrow}\lim\, \mathcal A_{X_i}$. as in \ref{4.1}, $\hat X={\rm Spec}\, \hat{\mathcal A}_X$. If $G$ is a $K$-scheme on groups, consider as a rational point the identity element.

Theorems \ref{universal} and \ref{4.4b} remain valid if $\bar X$ and $\bar G$ are substituted by $\hat X$ and $\hat G$. If, in addition, $C^*$ is assumed to be local, then the hypothesis of $X$ being compact and separated is no longer necessary  in \ref{8} and \ref{1.3}.

Finally, given an algebraic abelian group $G$, over an algebraically closed field $K$, it easy to prove that
$\bar G = \hat G\times G(K)$, where $G(K)$ is the constant functor $G(K)$,
that is, consider on $G(K)$ the discrete topology, then $G(K)(S):=\Aplic(\Spec S, G(K))$, for all commutative $K$-algebra.
\end{note}

\begin{proposition} \label{3.8} Let $\Spec \mathcal C^*$ be a formal monoid and $D_{\mathcal C}=\{w\in C^*\colon w(I)=0$ for some bilateral ideal $I\subset C$ of finite codimension$\}\subset C^*$.
Then, $\Spec D_{\mathcal C}$ is an affine monoid scheme and
$$\Hom_{mon}(\Spec \mathcal C^*, \Spec A)=
\Hom_{mon}( \Spec D_{\mathcal C},\Spec A)$$ for every affine monoid scheme $\Spec A$.
\end{proposition}

\begin{proof} Observe that $D_{\mathcal C}=\ilim{I} (C/I)^*$ and $\mathcal D_{\mathcal C}^*=\bar{\mathcal C}$. Then,
$$\aligned \Hom_{mon} & (\Spec \mathcal C^*, \Spec A) \overset{\text{\ref{5.3n}}}=\Hom_{bialg}(\mathcal A,\mathcal C^*)\,\overset{\text{\ref{n5.27}}}=\,\Hom_{bialg}(\mathcal C, {\mathcal A}^*)\\ &\overset{\text{\ref{n5.30}}} =\,
\Hom_{bialg}(\bar{\mathcal C},\mathcal A^*) \overset{\text{\ref{n5.27}}}=\Hom_{bialg}({\mathcal A},\mathcal D_{\mathcal C} ) \overset{\text{\ref{5.3n}}}=
\Hom_{mon}( \Spec D_{\mathcal C},\Spec A)\endaligned$$

\end{proof}

\begin{note}  Let $X$ be a $K$-scheme and $A_X$ the ring of functions of $X$.
The set $D_X$ of distributions of $X$ of finite support is said to be
$D_X:=\{w\in A^*_X\colon w$ factorizes through a finite quotient algebra of $A_X\}$. Obviously, $\mathcal D_X^*=\bar{\mathcal A}_X$ and $\Spec
\mathcal D_X^*=\bar X$.\end{note}

If $\Spec\mathcal C^*$ is an abelian formal monoid then $G=\Spec C$ is an affine abelian monoid scheme and $D_{\mathcal C}=D_G$, then

\begin{equation} \label{1.9} \Hom_{mon}(G^\vee, \Spec A)=
\Hom_{mon}( \Spec D_{G},\Spec A)\end{equation}  for every affine monoid scheme $\Spec A$.

\subsection{Functorial Cartier Duality}

\begin{definition} Let $\mathbb G$ be a functor of abelian monoids.
$\mathbb G^\vee := {\mathbb Hom}_{mon} (\mathbb G,{\mathcal R})$ (where we regard ${\mathcal R}$
as a monoid with its multiplication) is said to be the dual monoid of $\mathbb G$.
\end{definition}

If $\mathbb G$ is a functor of groups, then $\mathbb G^\vee = {\mathbb Hom}_{grp} (\mathbb G, G_m^\cdot)$.

\begin{theorem}\label{dual}
Assume that ${\mathbb  G}$ is a functor of abelian
monoids with a reflexive functor of functions. Then,
${\mathbb  G}^\vee={\rm Spec}\, (\mathbb  A_{\mathbb G}^*)$
(in particular, this equality shows that ${\rm Spec}\,\mathbb
A_{\mathbb G}^*$ is a functor of abelian monoids).
\end{theorem}

\begin{proof}
$  \mathbb G^\vee= {\mathbb Hom}_{mon} (
{\mathbb G},{\mathcal R})\overset{\text{\ref{2.5}}}={\mathbb Hom}_{\mathcal R-alg}
(\mathbb  A_{\mathbb G}^*,{\mathcal R})={\rm Spec}\, (\mathbb  A_{\mathbb G}^*).$
\end{proof}

\begin{note}\label{comentario}
Explicitly, ${\rm Spec}\, \mathbb A_{\mathbb G}^* = {\mathbb Hom}_{mon} ({\mathbb G},\mathcal R)$, $\phi \mapsto \tilde{\phi}$,
where $\tilde{\phi} (x)$ $ = \phi(x)$, for every $\phi \in {\rm
Spec}\, \mathbb A_{\mathbb G}^* = {\mathbb Hom}_{\mathcal R-alg} (\mathbb  A_{\mathbb G}^*,{\mathcal
R})$ and $x \in {\mathbb G} \to \mathbb A_{\mathbb G}^*$.

$\mathbb G^\vee$ is a functor of abelian
monoids ($(f \cdot f')(g) := f(g) \cdot f'(g)$, for every $f,f'\in
\mathbb G^\vee$ and $g \in \mathbb G$), the inclusion $\mathbb G^\vee = {\mathbb Hom}_{mon}
(\mathbb G,{\mathcal R}) \subset {\mathbb Hom} (\mathbb G,{\mathcal R})={\mathbb A_{\mathbb G}}$
is a morphism of monoids and the diagram $$\xymatrix{ \mathbb G^\vee \ar@{^{(}->}[r] \ar@{=}[d] & {\mathbb Hom}({\mathbb G}, {\mathcal R}) = \mathbb A_{\mathbb G}
\ar@{=}[d] \\ {\rm Spec}\, {\mathbb A_{\mathbb G}^*} \ar@{^{(}->}[r] & {\mathbb Hom}_{\mathcal R}(\mathbb A_{\mathbb G}^*,{\mathcal R})=\mathbb A_{\mathbb G}^{**}}$$ is commutative.
\end{note}

\begin{theorem} \label{Cartier}
The category of abelian affine $R$-monoid schemes  $G={\rm Spec}\, A$ is
anti-equiva\-lent to the category of abelian formal monoids ${\rm Spec}\, {\mathcal
A}^*$ (we assume the $R$-modules $A$ are
projective).
\end{theorem}

\begin{proof} The functors $\Spec A=\mathbb G\rightsquigarrow \mathbb G^\vee=\Spec\mathcal A^*$ and $\Spec\mathcal A^*=\mathbb G\rightsquigarrow \mathbb G^\vee=\Spec A$ establish the anti-equivalence between the category of abelian affine $R$-monoid schemes and the category of abelian formal monoids:

 The morphism $\mathbb G\overset{\delta}\to \mathbb G^{\vee\vee}$, $g \mapsto
\delta_g$, where $\delta_g (f) := f(g)$ for every $f \in \mathbb G^\vee $, is an
isomorphism: It is easy to check that the diagram
$$\xymatrix{{\rm Spec}\, {\mathbb A^{**}} \ar[r]^-\sim_-{\text{\ref{dual}}} & ({\rm
Spec} \, {\mathbb A^*})^\vee
& ({\rm Spec}\, {\mathbb A})^{\vee\vee}
\ar[l]_-\sim^-{\text{\ref{dual}}}\\ & {\rm Spec}\, {\mathbb A}
\ar@{=}[ul] \ar[ur]^-{\delta}& }$$ is commutative.

${\rm Hom}_{mon} (\mathbb G_1,\mathbb G_2) = {\rm Hom}_{mon}
(\mathbb G_2^\vee,\mathbb G_1^\vee)$:
Every morphism of monoids $\mathbb G_1 \to \mathbb  G_2$, taking ${\mathbb Hom}_{mon} (-,{\mathcal R})$, defines a morphism $\mathbb G_2^\vee \to \mathbb G_1^\vee$.
Taking ${\mathbb Hom}_{mon}(-,{\mathcal R})$ we get the original
morphism $\mathbb G_1 \to \mathbb G_2$, as it is easy to check.

\end{proof}

In \cite[Ch. I, \S 2, 14]{dieudonne}, it is given the Cartier Duality (formal schemes are certain functors over the category of commutative linearly compact algebras over a field).

Assume $G=\Spec A$ and $G'=\Spec B$ are commutative affine monoid schemes, then

\begin{equation} \label{1.92} \aligned
\Hom_{mon}( \Spec D_{G},G') & \overset{\text{Eq.\ref{1.9}}}=\Hom_{mon}(G^\vee, G') =
\Hom_{mon}(G'^\vee, G)\\ & \overset{\text{Eq.\ref{1.9}}}=\Hom_{mon}( \Spec D_{G'},G)\endaligned\end{equation}

\section{Examples of functors of monoids and Cartier duality}

\medskip
1.  {\bf Affine toric varieties.}
\medskip

Let $T$ be a set with structure of abelian (multiplicative)
monoid. Let $R$ be a field. The constant functor ${\bf T} = {\rm Spec}\, \prod_T
{\mathcal R}$ is an abelian formal monoid. The dual functor is the
abelian affine $R$-monoid scheme ${\bf T}^\vee = {\rm Spec}\, \oplus_T {\mathcal R}
= {\rm Spec}\, RT$.

An affine group scheme $G=\Spec A$ is linearly reductive if and only if $\mathcal A^*$ is linearly reductive. Since $\mathcal A^*$ is the inverse limit of its coherent algebra quotients
it is easy to check that $\mathcal A^*$ is linearly reductive if and only if it is a product of algebras of matrices. Then, $G=\Spec A$ is a linearly reductive commutative group scheme if and only if $\mathcal A^*=\prod_T \mathcal K$.

We will say that an abelian monoid $T$ is standard if it is finitely generated, its associated group $G$ is torsion-free and the natural morphism $T \to G$ is injective (in the literature, see \cite[6.1]{B}, it is called affine monoid). It is easy to prove that $T$ is standard if and only if $RT=\oplus_T R$ is a  finitely generated domain over $R$.

{\bf Theorem:}
The category of abelian monoids (resp. finitely generated monoids, standard monoids) is anti-equivalent to the category of affine  semisimple abelian monoid schemes (resp. algebraic affine semisimple abelian monoids, integral algebraic affine semisimple abelian monoids).

 If $T$ is standard then $G = \mathbb{Z}^n$ and the morphism $T \to G$ induces a morphism $G^n_m \to { \bf T}^\vee$. In particular, $G^n_m$ operates on ${\bf T}^\vee$. Furthermore, as $RG$ is the localization of $RT$ by the algebraically closed system $T$, the morphism $G^n_m \to {\bf T}^\vee$ is an open injection. We will say that an integral affine algebraic variety on which the torus operates with a dense orbit is an affine toric variety. It is easy to prove that there exists a one-to-one correspondence between affine toric varieties with a fixed point whose orbit is transitive and dense, and standard monoids.

\medskip
2.  {\bf Group $\mathbb Z$.}
\medskip

Obviously, $\mathbb Z^\vee=G_m$ and by Cartier duality $G_m^\vee=\mathbb Z$, which is a formal group.

\begin{proposition} Let $G$ be an affine  $K$-group scheme. Then,
$$\{\text{Rational points of}\, G\}=\Hom_{mon}(\mathbb Z,G) \overset{\text{Eq.\ref{1.9}}}=\Hom_{mon}(\Spec {D_{G_m}}, G)$$
If $G$ is an affine  commutative $K$-group scheme, by Equation \ref{1.92},
$$\{\text{Rational points of}\, G\}=\Hom_{mon}(\Spec {D_{G}}, G_m)$$\end{proposition}

3. {\bf Functor $\rad\mathcal K$:}
\medskip

Let $\rad\mathcal K\subset G_a$ be the covariant subfunctor of groups
defined by
$\rad{\mathcal K}(S):=\rad S=\{s\in S| s^n=0,\,\text{para algún }n\in\mathbb N\}$.

Let $K[x]/(x^n)\to K[x]/(x^m)$, $n\geq m$ be the natural quotient morphisms. Consider the projective system of $K$-álgebras $\{K[x]/(x^n),\, n\in \mathbb N\}$ and the inductive system of affine schemes $\{\Spec K[x]/(x^n),\, n\in \mathbb N\}$.

$\hat G_a=\ilim{n} \Spec{\mathcal K[x]/(x^n)}\overset{\text{\ref{4.4}}}=\Spec\plim{n}{\mathcal K[x]/(x^n)}=\Spec{\mathcal K[[x]]}$.

\begin{proposition} It holds that $\rad {\mathcal K} =\hat G_a=\Spec{\mathcal K[[x]]}$.\end{proposition}

\begin{proof} $(\Spec K[x]/(x^n))^\cdot(S)=\Hom_{K-alg}(K[x]/(x^n),
S)=\{s\in S| s^n=0\}$. Therefore,
$$\hat G_a(S)=\ilim{n} (\Spec
K[x]/(x^n))^\cdot(S)=\{s\in S| s^n=0,\,\text{para algún
}n\in\mathbb N\}=\rad{\mathcal K}(S)$$
\end{proof}

\begin{notation} Assume $\charc K=p$.
We denote $\alpha_{n}:=\Spec K[x]/(x^{p^n})\subset G_a$.
Given $\mu\in \alpha_1(S)\subset S$ (then $\mu^{p}=0$) we denote
$e^\mu:=\sum_{i=0}^{p-1} \mu^i/i!\in S$. Observe that $e^{\mu+\mu'}=e^\mu\cdot e^{\mu'}$.

We denote $\mu_n=\Spec K[x]/(x^{p^n}-1)\subset G_m$.
\end{notation}

Assume $\charc K=p$. Given $(\lambda_n)\in \prod_\mathbb N\alpha_1$
 the morphism $$ \rad{\mathcal K} \to  G_m,\,\, \mu \mapsto  e^{\lambda_0\mu}\cdot e^{\lambda_1\mu^p}\cdots e^{\lambda_n\cdot\mu^{p^n}}\cdots$$
 is a morphism of functors of groups: observe that the polynomial part of degree less than $p^{n+1}$ of the series $v(x):=e^{\lambda_0x}\cdot e^{\lambda_1x^p}\cdots e^{\lambda_n x^{p^n}}\cdots$ is
$e^{\lambda_0x}\cdot e^{\lambda_1x^p}\cdots e^{\lambda_n x^{p^n}}$.
The coefficient of  $x^{p^n}$ of $v(x)$ is $\lambda_n$, then $\prod_\mathbb N\alpha_1\subseteq (\rad{\mathcal K})^\vee$.

\begin{theorem} \label{1.16} \begin{enumerate} \item If $\charc K =0$, $ (\rad{\mathcal K})^\vee=G_a$.

\item If $\charc K=p>0$,
$(\rad{\mathcal K})^\vee=\prod_\mathbb N\alpha_1$,
as $K$-schemes.\end{enumerate}
\end{theorem}

\begin{proof} ${\mathbb Hom}(\rad{\mathcal K},{\mathcal K})=\mathcal K[[x]]$.
Hence,
$${\mathbb Hom}_{mon}(\rad{\mathcal K},{\mathcal K})=\{s(x)\in {\mathcal K[[x]]}\colon \begin{array}{l} s(\mu+\mu')=s(\mu)\cdot s(\mu'), \forall\, \mu,\mu'\in \hat G_a=\rad{\mathcal K}\\ s(0)=1\end{array}\}$$
In characteristic zero is easy to prove that $s(x)=e^{\lambda x}$, for some $\lambda\in \mathcal K$, then $ (\rad{\mathcal K})^\vee=G_a$.

Now, in characteristic $p>0$.
If the coefficient of $x$ of $s(x)$ is zero then is easy to prove that  $s(x)=t(x^p)$, with  $t(x)\in {\mathbb Hom}_{mon}(\rad{\mathcal K},{\mathcal K})$. If the coefficient of $x$ of $s(x)$ is $\lambda_0$ then the coefficient of
 $x$ of $(e^{\lambda_0x})^{-1}\cdot s(x)$ is zero, then
$(e^{\lambda_0x})^{-1}\cdot s(x)=t(x^p)$ and $s(x)=e^{\lambda_0x}\cdot t(x^p)$.
Likewise, $t(x)=e^{\lambda_1x}\cdot u(x^p)$ and $s(x)=e^{\lambda_0x}\cdot e^{\lambda_1x^p}
\cdot u(x^{p^2})$, etc. In conclusion, $\{(\lambda_n)\in \prod_\mathbb N\alpha_1\}= (\rad{\mathcal K})^\vee$

 \end{proof}

Assume $\charc K=p$.  $(\rad{\mathcal K})^\vee=\prod_\mathbb N\alpha_1=\Spec \mathcal K[x_0,\ldots,x_n,\ldots]/(x_0^p,\ldots,x_n^p,\ldots)$,
 as schemes. Observe that $$e^{\lambda_n\cdot\mu^{p^n}}\cdot e^{\lambda'_n\cdot\mu^{p^n}}=
e^{(\lambda_n+\lambda'_n)\cdot\mu^{p^n}}\cdot e^{\sum_{i=1}^{p-1}\frac{\lambda^i_n}{i!}\frac{{\lambda'}^{p-i}_n}{(p-i)!}\cdot\mu^{p^{n+1}}}.$$ Then, if $*$ denotes the operation of  $(\rad{\mathcal K})^\vee=\prod_\mathbb N\alpha_1$, then
$$(0,\cdots,0,\lambda_n,0,\cdots)*((0,\cdots,0,\lambda'_n,0,\cdots)=
(0,\cdots,0,\lambda_n+\lambda'_n,\sum_{i=1}^{p-1}\frac{\lambda^i_n}{i!}\frac{{\lambda'}^{p-i}_n}{(p-i)!},0,\cdots)$$

We have the natural inclusion $(\rad{\mathcal K})^\vee=\Spec \mathcal K[[x]]^*\hookrightarrow \mathcal K[[x]]$, $(\lambda_n)_{n\in\mathbb N}\mapsto
e^{\lambda_0x}\cdots e^{\lambda_nx^{p^n}}\cdots$. We also have the natural inclusion
$$\begin{array}{rcl} \rad{\mathcal K} & \hookrightarrow & {\mathbb Hom}((\rad{\mathcal K})^\vee,{\mathcal K})=\mathcal K[x_0,\ldots,x_n,\ldots]/(x_0^p,\ldots,x_n^p,\ldots)\\ \mu & \mapsto & e^{x_0\mu}\cdots e^{x_n\mu^{p^n}}\cdots\end{array} $$

\begin{note} Assume $\charc K=0$. $(\rad{\mathcal K})^\vee=G_a$ does not contain any proper subgroup. Then, the quotient groups of $(\rad{\mathcal K})^\vee$ are $G_a$ and the trivial group.

Assume $\charc K=p$. The finite subgroups of $G_a$ are $\{\alpha_n\}_{n\in\mathbb N}$ and
$\rad{\mathcal K}=\ilim{n} \alpha_n$. Then, $(\rad{\mathcal K})^\vee=\plim{n} \alpha_n^\vee$.
If $G$ is an algebraic group then
$$\Hom_{grp}((\rad{\mathcal K})^\vee, G)=
\Hom_{grp}(\plim{n} \alpha_n^\vee, G)=\ilim{n} \Hom_{grp}(\alpha_n^\vee,G)$$
That is, every morphism $(\rad{\mathcal K})^\vee\to G$ is the composition  of a projection
$(\rad{\mathcal K})^\vee\to\alpha_n^\vee$ and an injective morphism $\alpha_n^\vee\hookrightarrow G$.
\end{note}

\begin{theorem} In characteristic $p>0$, $(\alpha_n)^\vee=\prod_n\alpha_1$, as $K$-schemes. Specifically,
$$(\alpha_n)^\vee=\{ (\lambda_n)\in \prod^n\alpha_1\}=\left\{\begin{array}{rll} \alpha_n & \to & G_m\\ \mu & \mapsto & e^{\lambda_0\mu}\cdot e^{\lambda_1\mu^p}\cdots e^{\lambda_n\cdot\mu^{p^n}}\end{array}\right\}$$\end{theorem}

\begin{proof} Proceed as in Theorem \ref{1.16}.\end{proof}

 The subgroups of $\alpha_n$ are the obvious subgroups  $\alpha_r$, $r\leq n$, then, by Cartier duality the quotients of $\alpha_n^\vee$ are the obvious quotients $\alpha_r^\vee$, $r\leq n$.

\begin{proposition} Assume $\charc K=0$.  Then, $\rad\mathcal K=\hat G_m$.
\end{proposition}

\begin{proof}  The inverse morphism of the morphism $\rad\mathcal K\to \hat G_m,\,\mu\mapsto e^\mu$ is the morphism $\hat G_m\to \rad\mathcal K$, $\mu\mapsto \ln (1+(\mu-1)):=\sum_{i>0} (-1)^{i+1} (\mu-1)^i/i!$.\end{proof}


In characteristic $p>0$, $\hat G_m = \ilim{n} \mu_{n}$ as functor of groups.

\medskip
4. {\bf Functor of monoids $G_a^\vee$:}
\medskip

 Let $G_a={\rm Spec}\, K[x]$ be the additive group.

\begin{theorem} Assume ${\rm char}\, K=0$. Then
$$G_a^\vee=\rad{\mathcal K}$$\end{theorem}

\begin{proof}  We know that $(\rad{\mathcal K})^\vee=G_a$, by \ref{1.16}.
By Cartier duality, $G_a^\vee=\rad{\mathcal K}$.\end{proof}

Explicitly, $\mu\in\rad{\mathcal K}$ defines the morphism $G_a\to G_m$,
$\lambda\mapsto e^{\lambda \mu}$. Then, we have the natural inclusion
$\rad{\mathcal K}\hookrightarrow \mathcal K[x]$, $\mu\mapsto e^{\mu x}$.

\begin{theorem} \label{7.9} Assume $\charc K =p>0$. Then, $G_a^\vee$ is isomorphic to  $\oplus_\mathbb N \alpha_1$ as  functors of sets.  Specifically,
$$G_a^\vee=\{(\lambda_n)\in \oplus_\mathbb N \alpha_1\}=
\left\{\begin{array}{rll} G_a & \to  & G_m\\ \alpha & \mapsto &
e^{\lambda_0\alpha}\cdot e^{\lambda_1\alpha^p}\cdots e^{\lambda_n\cdot \alpha^{p^n}}
\end{array}\right.$$
\end{theorem}

\begin{proof} Proceed as in Theorem \ref{1.16}.\end{proof}

Write $\mathcal K[[x_0,\ldots,x_n,\ldots]]/(x_0^p,\ldots,x_n^p,\ldots):=\plim{n}
\mathcal K[x_0,\ldots,x_n]/(x_0^p,\ldots,x_n^p)$, then
$G_a^\vee=\oplus_\mathbb N \alpha_1=\Spec  \mathcal K[[x_0,\ldots,x_n,\ldots]]/(x_0^p,\ldots,x_n^p,\ldots)$, as functors of sets. The dual morphism of natural the inclusion $\rad \mathcal K\subset G_a$ is the obvious morphism $G_a^\vee=\oplus_\mathbb N\alpha_1\subset \prod_\mathbb N\alpha_1=(\rad\mathcal K)^\vee$.

We have the natural inclusion $G_a^\vee\hookrightarrow \mathcal K[x]$, $(\lambda_0,\cdots,\lambda_n)\mapsto
e^{\lambda_0x}\cdots e^{\lambda_nx^{p^n}}$.

\begin{note} Let us now assume  that $\charc K=0$. Let $\mathbb Q$ be the
constant functor of groups    $\mathbb Q$ (with the addition
operation), that is, denote $\mathbb Q$ the discrete topological
space, then  $\mathbb Q(S):=\Aplic(\Spec S,\mathbb Q)$, for all
commutative $K$-algebras $S$.  Let $G_M=\mathbb Q^\vee$ be the
dual group of $\mathbb Q$, then $G_M=\Spec (K[\mathbb
Q]=K[e^{rx}]_{r\in\mathbb Q})$. $\mathbb Q$ is the direct limit of
its finite $\mathbb Z$-submodules, $\mathbb Q=\ilim{r\in \mathbb
Q^+} \mathbb Z\cdot r$, which are isomorphic to $\mathbb Z$, then
$G_M=\plim{r\in \mathbb Q^+} \Spec K[e^{rx},e^{-rx}]$ and $\Spec
K[e^{rx},e^{-rx}]\simeq G_m$. Giving a point $\alpha\in G_M$ is
equal to giving a point $\alpha\in G_m$, and to determining
$\alpha^r$,for all $r\in\mathbb Q$.  If $G$ is an algebraic group,
then
$$\Hom_{mon}(G_M,G)  =\ilim{r\in\mathbb Q^+}\Hom_{mon}(\Spec
K[e^{rx},e^{-rx}], G)$$
That is, every morphism from $G_M$ to an algebraic group factors through a projection
$G_M\to G_m$, $\alpha\mapsto \alpha^r$, for some $r\in\mathbb Q$.

If $E=\oplus_J\mathbb Q$ is a $\mathbb Q$ is a vector space, let us also write  $E$ for the
constant functor of monoids $E$, then
$E^\vee:={\mathbb Hom}_{mon}( E,G_m)= \prod_J G_M$. Every morphism of $ E^\vee$ on an algebraic group factors through $G_m^n$, for some $n$.

Let us now assume  that $\charc K=p>0$. If $E=\oplus_J\mathbb
Z/p\mathbb Z$ is a $\mathbb Z/p\mathbb Z$-vector space then
$E^\vee:={\mathbb Hom}_{mon}(E,G_m)= \prod_J \mu_1$. Every
morphism from $E^\vee$ to an algebraic group factors through
$\mu_1^n$, for some $n\in \mathbb N$.\end{note}

\begin{theorem} \label{37} Let $G_a=\Spec K[x]$ and let $D_{G_a}$
be the set of distributions with finite support of $G_a$.
Let $G_a(K)$ be the constant functor of groups, $K$.
Assume $K$ is an algebraically closed field. Then,
\begin{enumerate}

\item $\bar G_a=\hat G_a\times G_a(K)=\rad\mathcal K\times G_a(K)$ and $\Spec D_{G_a} = (\rad\mathcal K)^\vee\times G_a(K)^\vee$.

\item
It holds that
$$\Hom_{mon}(G_a^\vee, G)\overset{\text{Eq.\ref{1.9}}}=\Hom_{mon}((\rad\mathcal K)^\vee \times
G_a(K)^{\vee},G),$$ for all affine group scheme  $G$.

\item The scheme-theoretic image of every morphism from $G_a^\vee$ to an affine algebraic group is isomorphic to

$$\begin{array}{ll} G_a^\delta\times G_m^n,\,
\delta=0,1, & \text{ if }\charc K=0 \\ \alpha_r^\vee\times \mu_1^s, & \text{ if } \charc K=p>0\end{array}$$

\end{enumerate} \end{theorem}

\section{Exponential map in arbitrary characteristic}

\subsection{Vector fields} ${}$

\medskip

Let $K[\epsilon]=K[x]/(x^2)$. Given a functor of sets  $\mathbb X$
and the morphism of $K$-algebras $K[\epsilon]\to K$,
$\epsilon\mapsto 0$, let $0^*\colon \mathbb X(K[\epsilon])\to
\mathbb X(K)$ be the induced morphism. Given $x\in \mathbb X(K)$,
let
$$T_x \mathbb X:=\{y\in \mathbb X(K[\epsilon]), \text{ such
that } 0^*y=x\}$$

\begin{proposition} Let $X$ be a $K$-scheme and $x\in X$ a rational point.
Then, $T_x X= T_x \bar X$.\end{proposition}

\begin{proof} Observe that $\bar X(K[\epsilon])=X(K[\epsilon])$ y
$\bar X(K)=X(K)$.

\end{proof}
More generally, given a commutative $K$-algebra $S$ and the morphism of $S$-algebras
$S[\epsilon]=S[x]/(x^2)\to S$,
$\epsilon\mapsto 0$ let $0^*\colon \mathbb X(S[\epsilon])\to \mathbb X(S)$ be the induced morphism. Given $x\in \mathbb X(S)$, let
$$T_x \mathbb X:=\{y\in \mathbb X(S[\epsilon]), \text{ such
that } 0^*y=x\}$$ be the vector space of tangent vectors of $\mathbb X$ at the point $x$.

Assume $\mathbb X=\Spec \mathbb A$. Then $\mathbb X(S)=\Hom(\Spec S,\Spec \mathbb A)=
\Hom_{\mathcal K-alg}(\mathbb A,\mathcal S)$ and $\mathbb X(S[\epsilon])=\Hom_{\mathcal K-alg}(\mathbb A,\mathcal S[\epsilon])$. Then, given $x\in \mathbb X(S)=\Hom_{\mathcal K-alg}(\mathbb A,\mathcal S)$,  by the standard arguments
$$T_x\Spec\mathbb A=\Der_{\mathcal K}(\mathbb A, \mathcal S)$$

Given a morphism of functors of sets $\phi\colon \mathbb X\to \mathbb Y$, $x\in \mathbb X(S)$ and $D_x\in T_x\mathbb X\subset \mathbb X(S[\epsilon])$
then $\phi(D_x)\in T_{\phi(x)}\mathbb Y\subset \mathbb Y(S[\epsilon])$.

\begin{definition} We will say that $T\mathbb X:=\mathbb Hom(\Spec K[\epsilon],\mathbb X)$
is the tangent bundle of $\mathbb X$. We have a natural morphism
$T\mathbb X\overset{0^*}\to \mathbb X$, $D_x\mapsto 0^*D_x=x$. We
will say that $\Der \mathbb X:=\Hom_{\mathbb X}(\mathbb X,
T\mathbb X)$ is the set of derivations of $\mathbb X$ (or the set
of vector fields of $\mathbb X$).\end{definition}

Giving $D\in \Der \mathbb X$ is equivalent to giving for every
point  $x\in \mathbb X(S)$ (and every $S$) a tangent vector $D_x\in T_x \mathbb X\subset
\mathbb X(S[\epsilon])$, functorially.

Let $G_a=\mathcal K$ be the additive group, then $$TG_a=\mathbb Hom(\Spec K[\epsilon],G_a)=
\mathbb Hom_{\mathcal K-alg}(\mathcal K[x],\mathcal K[\epsilon])=\mathcal K[\epsilon]$$

\begin{proposition} Let $\mathbb X$ be an affine functor, then it holds that
$$\Der \mathbb X=\Der_{\mathcal K}(\mathbb A_{\mathbb X},\mathbb A_{\mathbb X})$$
\end{proposition}

\begin{proof} First, observe that
$$\aligned \mathbb A_{\mathbb X\times\Spec K[\epsilon]} & =
\mathbb Hom(\mathbb X\times \Spec K[\epsilon],\mathcal K)=
\mathbb Hom(\mathbb X, \mathbb Hom(\Spec K[\epsilon],\mathcal K))
\\ & =\mathbb Hom(\mathbb X, \mathcal K[\epsilon])=\mathbb A_{\mathbb X}[\epsilon]\endaligned$$
Now,
$$\aligned \Der \mathbb X & =  \Hom_{\mathbb X}(\mathbb X, T\mathbb X)
=\Hom_{\mathbb X}(\mathbb X, \mathbb Hom(\Spec K[\epsilon]
,\mathbb X))\\ & =\{f\in \Hom(\mathbb X\times \Spec K[\epsilon]
,\mathbb X)\colon f(x,0)=x,\,\forall x\in\mathbb X\}\\ & = \{h\in
\Hom_{\mathcal K-alg}(\mathbb A_{\mathbb X},\mathbb A_{\mathbb
X}[\epsilon])\colon h=\Id\mod(\epsilon)\}= \Der_{\mathcal
K}(\mathbb A_{\mathbb X},\mathbb A_{\mathbb X})\endaligned$$

\end{proof}

Given a $K$-scheme $X$, it is well known that $$\Der X = T_{\Id}
{\mathbb End}\, X\subset End_{K[\epsilon]-sch}(X\times_K \Spec
{K[\epsilon]}),\quad D\mapsto e^{\epsilon D},$$ (topologically
$e^{\epsilon D}$ is the identity morphism and over the ring of
functions $e^{\epsilon D}(a+b\epsilon):=a+\epsilon \cdot
D(a)+b\epsilon$, for all $a+b\epsilon\in \mathcal O_X[\epsilon]$).

\begin{theorem} Let $\mathbb X$ be a functor of sets. It holds that
$$T_{\Id} {\mathbb End}\, \mathbb X=\Der \mathbb X$$\end{theorem}

\begin{proof}  It is a consequence of the equalities
$$\aligned \Hom(\mathbb X_{|K[\epsilon]},\mathbb X_{|K[\epsilon]}) & =
\Hom_{\Spec K[\epsilon]}(\mathbb X\times \Spec K[\epsilon],\mathbb X\times \Spec K[\epsilon])\\ & = \Hom(\mathbb X\times \Spec K[\epsilon],\mathbb X) =
\Hom(\mathbb X,\mathbb Hom(\Spec K[\epsilon],\mathbb X))\endaligned$$

\end{proof}

Giving $D\in \Der\mathbb X$, we denote $\mathbb X\to  \mathbb
Hom(\Spec K[\epsilon],\mathbb X)$, $x\mapsto D_x$ the
corresponding morphism and we   denote $e^{\epsilon D}\colon
\mathbb X_{|K[\epsilon]}\to \mathbb X_{|K[\epsilon]}$ the
corresponding endomorphism.
 Given $x\in \mathbb X(S)\subset \mathbb X(S[\epsilon])$, it is easy to check that $e^{\epsilon D}(x)=D_x$.

\begin{definition} Let $\mathbb G$ be a functor of groups. We say that a vector field $D\in\Der \mathbb G$ (for each point $g\in \mathbb G(S)$
we have a tangent vector $D_g$ at $g$)  is
$\mathbb G$-invariant if $g'\cdot D_{g}=D_{g'\cdot g}$ for all pair of points
$g,g'$ of $\mathbb G$. We will say that $\Derinv \mathbb G$ is the set of all invariant
vector fields of $\mathbb G$.\end{definition}

\begin{proposition} Let  $\mathbb G$ be a functor of groups and let $e\in \mathbb G(K)$ be the
identity element. It holds that $\Derinv \mathbb G = T_e\mathbb G$.\end{proposition}

\begin{proof} The morphisms $\Derinv \mathbb G \to T_e\mathbb G$, $D\mapsto D_e$, $T_e\mathbb G\to
\Derinv \mathbb G$, $D_e\mapsto D$, where $D_g:=g\cdot D_e$ for all
$g\in \mathbb G$, are mutually inverses.\end{proof}

\begin{definition} Given $D_e\in T_e\mathbb G$, we will say that $D\in\Derinv\mathbb G$, such that
$D_g:=g\cdot D_e$, for all $g\in\mathbb G$, is the invariant field associated with $D_e$.\end{definition}

\begin{proposition}  Let $f\colon \mathbb G_1\to \mathbb G_2$ be a morphism of functors of groups, $D_e\in T_e\mathbb G_1$
and $D'_e=f(D_e)$. If $D$ y $D'$ are the invariant fields associated with  $D_e$ and $D'_e$ respectively,
then  $f(D_g)=D'_{f(g)}$, for all $g\in \mathbb G_1$.\end{proposition}

\begin{proof}
$$f(D_g)=f(g\cdot D_e)=f(g)\cdot
f(D_e)=f(g)\cdot (D'_e)=D'_{f(g)}$$ \end{proof}

\subsection{Analytic one-parameter group} ${}$

\medskip

\begin{theorem} \label{6.3} Let $\mathbb X$ be an affine functor of sets. It holds that
$$\Der(\mathbb X)=\Der_{\mathcal K}(\mathbb A_{\mathbb X},\mathbb A_{\mathbb X})=
\Hom_{mon}(G_a^\vee,{\mathbb End}\,\mathbb X)$$\end{theorem}

\begin{proof} Let $\mathbb M$ be a dual functor of $\mathcal K$-modules. By Theorem \ref{2.5}, giving a $\mathcal K[x]$-module structure
on ${\mathbb M}$ (equivalently, giving a linear endomorphism on
${\mathbb M}$) is equivalent to giving  a
$G_a^\vee$-module structure on $\mathbb M$. If ${\mathbb M}$ is a $\mathcal K[x]$-module then through the inclusion $G_a^\vee\hookrightarrow \mathcal \mathcal K[x]$,
$G_a^\vee$ operates on $\mathbb M$. Let us follow the notations
$$\mathbb End_K{\mathbb M}=\Hom_{K-alg}(\mathcal K[x],\mathbb End_K{\mathbb M})=\Hom_{mon}(G_a^\vee,{\mathbb End}_{\mathcal K}\mathbb M),\quad T\mapsto e^{xT}$$
 where $e^{xT}\colon G_a^\vee\to {\mathbb End}_{\mathcal K}{\mathbb M}$ is defined by  $e^{xT}(\mu)=e^{\mu T}$ (if $\charc K=0$, $\mu\in \rad\mathcal K=G_a^\vee$ and $e^{\mu T}:=\sum_i (\mu\cdot T)^i/i!$; if
$\charc K=p>0$, $\mu=(\mu_0,\ldots,\mu_r)\in \oplus_{\mathbb N}\alpha_1=G_a^\vee$ and
$e^{\mu T}:=e^{\mu_0 T}\cdot e^{\mu_1 T^p}\cdots e^{\mu_r T^{p^r}}$).

Likewise, giving a $\mathcal K[x_1,\ldots,x_n]$-structure on $\mathbb M$ is equivalent to
giving a $G_a^\vee\times\overset n\cdots\times G_a^\vee$-module structure on $\mathbb M$.

Let us follow the notations
$$\aligned (\mathbb End_K \mathbb M)^n\supseteq \Hom_{K-alg}(K[x_1,\ldots,x_n],\mathbb End_K\mathbb M)
 & = \Hom_{mon}((G_a^\vee)^n, {\mathbb End}_K \mathbb M),\quad
\\ (T_1,\ldots,T_n) & \overset{Not.} \mapsto e^{x_1T_1}\cdots e^{x_nT_n}\endaligned$$

If $T_1,T_2\in \mathbb End_K\mathbb M$ commute then
$e^{x(T_1+T_2)}=e^{xT_1}\cdot e^{xT_2}$.
Let  $D\in End_K\mathbb A_{\mathbb X}$, that is, $\mathbb A_{\mathbb X}$ is a $\mathcal K[x]$-module (hence it is a
$G_a^\vee$-module: given $\mu\in G_a^\vee$, $\mu\cdot a=e^{\mu D} a$).
Consider $D\otimes 1+1\otimes D\in \mathbb End_K (\mathbb A_{\mathbb X}\otimes \mathbb A_{\mathbb X})$, then
$\mathbb A_{\mathbb X}\otimes \mathbb A_{\mathbb X}$ is a $\mathcal K[x]$-module (hence it is a $G_a^\vee$-module: $\mu\cdot
(a\otimes b)= e^{\mu(D\otimes 1+1\otimes D)}(a\otimes b)=e^{\mu D}
a\otimes e^{\mu D} b$). Now , $D$ is a derivation if and only if
the morphism $\mathbb A_{\mathbb X}\otimes \mathbb A_{\mathbb X}\to \mathbb A_{\mathbb X}$, $a\otimes b\mapsto ab$ is a morphism of $\mathcal K[x]$-modules,
which is equivalent to say that the morphism is a morphism of
 $G_a^\vee$-módulos, that is, $G_a^\vee$ operates on $\mathbb A_{\mathbb X}$ by morphisms of $K$-algebras. That is,
the diagram

$$\xymatrix{ \mathbb End_K\mathbb A_{\mathbb X} \ar@{=}[r] & \Hom_{mon}(G_a^\vee,{\mathbb End}_K{\mathbb A_{\mathbb X}})\\ \Der_K \mathbb A_{\mathbb X} \ar@{=}[r] \ar@{^{(}->}[u] &
\Hom_{mon}(G_a^\vee,{\mathbb End}_{K-alg}{\mathbb  A_{\mathbb X}}) \ar@{^{(}->}[u]
}$$ is commutative.\end{proof}

\begin{theorem} \label{55} Let $X$ be a $K$-scheme. It holds that
$$\Hom_{mon}(G_a^\vee,{\mathbb End}\, X)=\Der X$$
\end{theorem}

\begin{proof} Given a point $\mu\in G_a^\vee(S)$ there exists a finite, local and rational $K$-subalgebra  $C\subset S$, such that $\mu\in G_a^\vee(C)$. Moreover, if $C\to K$ is the quotient morphism then $G_a^\vee(C)\to G_a^\vee(K)=\{0\}$, $\mu\mapsto 0$. Therefore, given a morphism $f\colon
G_a^\vee\to {\mathbb End} X$, then  $f(\mu)\in {\mathbb End} X(C)$, that is, $f(\mu)$ is an endomorphism $X_C\to X_C$.
Moreover, by base change $C\to K$, $f(\mu)=\Id$. Topologically $X_C=X$, then $f(\mu)$ is topologically equal to the identity morphism, hence it is affine.
 Let $\{U_i\}$ be a covering of $X$ by open affine sets. Giving a morphism $G_a^\vee\to {\mathbb End} X$
 is equivalent to giving morphisms $G_a^\vee\to {\mathbb End}\, U_i$ which coincide on $U_i\cap U_j$. That is, giving a morphism $G_a^\vee\to {\mathbb End} X$ is equivalent to giving derivations
on each $U_i$ which coincide on $U_i\cap U_j$. That is,
$$\Hom_{mon}(G_a^\vee,{\mathbb End} X)=\Der X$$

\end{proof}

Let $D$ be a vector field on an affine functor of sets $\mathbb X$.  Consider the morphism

$$exp_D\colon G_a^\vee \to {\mathbb End}\, \mathbb X, \quad \mu \mapsto e^{\mu D}$$
If $X$ is a $K$-scheme, then $e^{\mu D}$ is an affine morphism
because it is topologically the identity morphism. $e^{\mu D}\colon \mathbb X_{|C}\to \mathbb X_{|C}$ ($\mu\in G_a^\vee(C)$)   induces the morphism $e^{\mu D}\colon (\mathbb A_{\mathbb X})_{|C}
 \to (\mathbb A_{\mathbb X})_{|C}$ defined by
$$e^{\mu D}(a):=\left\{\begin{array}{ll} \sum_n \frac{\mu^n D^n(a)}{n!},\, & \text{ if
}\charc K=0\\ e^{\mu_0 D}(a)\cdots e^{\mu_n D^{p^n}}(a), & \text{ if}
\charc K=p>0\end{array}\right.$$

 If $X$ is a projective $K$-variety, then there exists a $K$-scheme whose functor of points is $\mathbb End\,X^\cdot$
 (see \cite{G}), and we will also denote it $\mathbb End\, X$. Thus, $exp_D\colon G_a^\vee\to\mathbb End\,X$ factors through a unique morphism $exp_D\colon \Spec K[x]^*\to\mathbb End\,X$.

\begin{definition} The morphism  $exp_D\colon G_a^\vee\to {\mathbb End}\, \mathbb X$ is said to be the analytic one-parameter group associated with $D$.

Let $y\in \mathbb X$ be a rational point. We will say that $exp_{D,y}\colon G_a^\vee\to \mathbb X$, $exp_{D,y}(\mu):=exp_D(\mu)(y)$ is the analytic integral curve of $D$ passing through $y$.\end{definition}

The induced morphism between the rings of functions by
$exp_{D,y}$ is $\mathbb A_{\mathbb X}\to\mathcal K[x]^*$,
$$a\mapsto \left\{\begin{array}{ll} \sum_i D^i(a)(y)/i!\cdot x^i\in \mathcal K[[x]],& \text{if }\charc K=0\\
(\sum_{i=0}^{p-1} D^i(a)(y)/i!\cdot x_0^i)\cdots (\sum_{i=0}^{p-1}
D^{ip^n}(a)(y)/i!\cdot x_n^i)\cdots, & \text{if }\charc
K=p>0\end{array}\right.
$$

If $X$ is a $K$-scheme, the morphism $exp_{D,y}\colon G_a^\vee\to X$ factors through $\Spec K[x]^*$ and we also write it
$exp_{D,y}\colon \Spec K[x]^*\to X$. If $U=\Spec A$ is an affine open set containing $y$, then the morphism induced by $exp_{D,y}$ between the rings of functions is written as it has been written above.

\begin{theorem} \label{8.4} Let $\mathbb G$ be an affine functor of groups or a $K$-group scheme. It holds that
$$\Hom_{grp}(G_a^\vee, \mathbb G)=\Derinv(\mathbb G)=T_{e} \mathbb G$$\end{theorem}

\begin{proof} 1. Given $g\in \mathbb G$, let $L_g\colon \mathbb G\to \mathbb G$
 be the morphism defined by $L_g(g')=gg'$. $\mathbb G$ operates on
 $\mathbb End\, \mathbb G$ by $g*f:=L_g\circ f\circ L_{g^{-1}}$, for all
 $g\in \mathbb G$ and $f\in \mathbb End\, \mathbb G$.
 The morphism $f$ is $\mathbb G$-invariant if and only if $f(g)=g\cdot f(e)$, for all $g\in \mathbb G$.
 Hence, $(\mathbb End\, \mathbb G)^{\mathbb G}=\mathbb G$.

2. $\mathbb G$ operates on $\Der \mathbb G$ by $(g*D)_h:=g\cdot D_{g^{-1}h}$
 for all $g,h\in \mathbb G$ and $D\in \Der \mathbb G$. Obviously,
 $(\Der \mathbb G)^{\mathbb G}=\Derinv \mathbb G$.

3.  $\mathbb G$ operates on $\Hom_{mon}(G_a^\vee, \mathbb End\,\mathbb G)$
by $(g*F)(\mu)=g*(F(\mu))$, for all $g\in\mathbb G$, $F\in \Hom_{mon}(G_a^\vee, \mathbb End\,\mathbb G)$ and $\mu \in G_a^\vee$.

4. Taking $\mathbb G$-invariants on $\Der(\mathbb G)=\Hom_{mon}(G_a^\vee,{\mathbb End}\,\mathbb G)$,
by 1.,2. and 3. we obtain
$$\Derinv(\mathbb G)=\Hom_{grp}(G_a^\vee, \mathbb G)$$

\end{proof}

\subsection{Existence and uniqueness of ``analytic'' solutions of an algebraic
differential equation in arbitrary characteristic} ${}$

\medskip

\begin{proposition} If $\mu\in G_a^\vee(S)$, then $$T_\mu G_a^\vee =\left\{\begin{array}{l} S, \text{ if } \charc K=0\\ \oplus_{\mathbb N} S,\text{ if }\charc
K=p>0\end{array}\right.$$ \end{proposition}

\begin{proof} $G_a^\vee(S[\epsilon])=\{\mu+\lambda\cdot \epsilon$, with $\mu\in G_a^\vee(S)$ and $\lambda\in S$ if $\charc K=0$,
$\lambda\in \oplus_\mathbb N S$ if $\charc K=p>0\}$.

\end{proof}

In characteristic zero, let us denote $\mu+\lambda\cdot \epsilon =
\lambda({\partial_{x_0}})_\mu$. In positive characteristic $p>0$,
$G_a^\vee=\oplus_\mathbb N \alpha_1=\Spec
K[[x_0,\ldots,x_n,\ldots]]/(x_0^p,\ldots,x_n^p,\ldots)$. Let us
denote $\mu+\lambda\cdot \epsilon
=\sum_i\lambda_i({\partial_{x_i}})_{\mu_i}$. The morphism
$G_a^\vee\to \Spec K[x]^*$ induces a morphism $T_\mu G_a^\vee\to
T_\mu \Spec K[x]^*$, that maps $\mu+\lambda\cdot \epsilon$ to
$\sum_i \lambda_i({\partial_{x_i}})_{\mu_i}$.

The linear map induced between the tangent spaces by $exp_D$ is:
$$T_0 G_a^\vee \overset{exp_{D}}\to T_{Id}
{\mathbb Aut } \mathbb X= \Der \mathbb X, \quad \epsilon =({\partial_{x_0}} )_0\overset{exp_{D}}\mapsto e^{\epsilon D}=D$$
Then, by Theorem \ref{55} we have the following theorem.

\begin{theorem}[existence and uniqueness\,] Let $\mathbb X$ be an affine functor of sets or a $K$-scheme. If $D$ is a vector field on $\mathbb X$, then $exp_{D
}$ is the only morphism of functors of groups  $f\colon G_a^\vee
\to {\mathbb Aut}\, \mathbb X$ such that
$f((\partial_{x_0})_0)=D$.\end{theorem}

Let $x\in\mathbb X(K)$ be a rational point. The morphism
$\chi_x\colon \mathbb End\,\mathbb X\to \mathbb X$,
$\chi_x(\tau)=\tau(x)$, maps $e^{\epsilon D}=D$ to $e^{\epsilon
D}(x)=D_x$. Then,
$$exp_{D,x}((\partial_{x_0})_0)=\chi_x(exp_D((\partial_{x_0})_0))=\chi_x(D)=D_x$$

\begin{notation} Let $\delta\in\Derinv G_a^\vee$ be the invariant field associated with
$(\partial_{x_0})_0$. That is,
$\delta_\mu=\mu*(\partial_{x_0})_0$, where $\mu\in G_a^\vee$ and
$*$ is the operation of $G_a^\vee$. Specifically,
$\delta=\partial_{x_0}+x_0^{p-1}/(p-1)
!\partial_{x_1}+x_0^{p-1}/(p-1)!x_1^{p-1}/(p-1)!\partial_{x_2}+\cdots$.\end{notation}

It can be checked that:
$\delta^{p^n}=\partial_{x_n}+x_n^{p-1}/(p-1)!\partial_{x_{n+1}}+
x_n^{p-1}/(p-1)!x_{n+1}^{p-1}/(p-1)!\partial_{x_{n+2}}+\cdots$. Therefore: ${(\delta^{p^n}})_0=(\partial_{x_n})_{0}$
and:
$$\Derinv G_a^\vee=\langle \delta,\ldots, \delta^{p^n},\ldots \rangle$$

Let $\mu\in G_a^\vee(S)$, $x\in\mathbb  X(S)$ and denote $*$ the operation of $G_a^\vee$. We have the commutative diagram
$$\xymatrix{ G_a^\vee \ar[r]^-{exp_{D,x}} \ar[d]^-{\mu*} & \mathbb X \ar[d]^-{exp_D(\mu)}\\ G_a^\vee \ar[r]^-{exp_{D,x}} & \mathbb X  } \xymatrix{ \delta_{0} \ar@{|->}[r] \ar@{|->}[d] & D_x \ar@{|->}[d]  \\  \delta_{\mu} \ar@{|->}[r] & exp_{D,x}(\delta_{\mu})= exp_D(\mu)(D_{x})\,\,\,\,\,\,\,}$$

Let $\epsilon=\delta_0\in T_0G_a^\vee$. From the commutative diagram

$$\xymatrix{(G_a^\vee)_{|K[\epsilon]} \ar[r]^-{exp_{D,x}} \ar[d]^-{\epsilon*} & \mathbb X_{|K[\epsilon]} \ar[d]^-{exp_D(\epsilon)=e^{\epsilon D}}\\ (G_a^\vee)_{|K[\epsilon]} \ar[r]^-{exp_{D,x}} & \mathbb X_{|K[\epsilon]}}\quad\xymatrix{\mu \ar@{|->}[r] \ar@{|->}[d] &
\,\,\,\,\,\,\,\,\, exp_{D,x}(\mu) \ar@{|->}[d]\\ \delta_\mu \ar@{|->}[r] & exp_{D,x}(\delta_\mu)=
D_{\exp_{D,x}(\mu)}\,\,\,\,\,\,\,\,\,\,\,\,\,\,\,\,\,\,\,\,}$$
we obtain that $exp_D(\mu)(D_x)=exp_{D,x}(\delta_{\mu})=D_{exp_{D,x}(\mu)}=D_{exp_D(\mu)(x)}$.

\begin{lemma} Let $a\in \mathbb A_{G_a^\vee}=\mathcal K[x]^*$. Then, $a\in \mathcal K$ $\iff$ $\delta(a)=0$ $\iff$
$\delta(a)(\mu)=0$ for all $\mu\in G_a^\vee$. \end{lemma}

\begin{proof} Obviously $a\in \mathbb Hom(\mathbb X,\mathcal K)=\mathbb A_{\mathbb X}$ is zero if and only if
$a(x)=0$ for all $x\in\mathbb X$.

Let us only prove that $a\in \mathcal K$ if $\delta(a)=0$, and
$\charc K=p$: Obviously $\delta^{p^n}(a)=0$ for all $n\in\mathbb
N$. As
$$\mathcal K[x]^*=\mathcal K[[x_0,\ldots,x_n,\ldots]]/(x_0^p,\ldots,x_n^p,\ldots)=
\{\sum_{\alpha\in\oplus_\mathbb N\mathbb Z/p\mathbb Z}\lambda_\alpha x^\alpha,\,\lambda_\alpha\in \mathcal K\}$$
and $T_\mu G_a^\vee=\langle \delta_\mu,
\ldots,(\delta^{p^n})_\mu,\ldots\rangle=\langle (\partial_{x_i})_\mu\rangle$, then $\partial_{x_i}a=0$, for all
$i\in\mathbb N$, and $a\in \mathcal K$.

\end{proof}

\begin{theorem}[existence and uniqueness\,] \label{63} Let $\mathbb X$ be an affine functor of sets or a $K$-scheme, let $x\in \mathbb X(K)$ be a rational point and let $D$ be a vector field on $\mathbb X$. Then,
$exp_{D,x}$ is the only morphism $f\colon G_a^\vee \to \mathbb X$
such that $f(0)=x$ and $f({\delta}_{\mu})=D_{f(\mu)}$, for every
point $\mu\in G_a^\vee$.\end{theorem}

\begin{proof} We already know that $exp_{D,x}(0)= e^{0D}(x)=\Id(x)=x$ and $exp_{D,x}(\delta_\mu)= D_{exp_{D,x}(\mu)}$.

We still have to prove the uniqueness of $f$.

\noindent Let $\tilde D$ be the only invariant field on ${\mathbb
Aut} X$ such that $\tilde D_{\Id}=e^{\epsilon D}$. Recall that the
morphism $\chi_x\colon \mathbb Aut X\to X$, $\chi_x(\phi)=\phi(x)$
maps $\tilde D_{\Id}=e^{\epsilon D}$ to $e^{\epsilon D}(x)=D_x$.
We have the commutative diagram
$$\xymatrix{\mathbb Aut\,\mathbb X \ar[r]^-{\chi_x} \ar[d]^-{\phi\cdot} & \mathbb X\ar[d]^-\phi
\\ \mathbb Aut\,\mathbb X \ar[r]^-{\chi_x} & \mathbb X}\quad \xymatrix{\tilde D_{\Id} \ar@{|->}[r]^-{\chi_x} \ar@{|->}[d]^-{\phi\cdot} & D_x  \ar@{|->}[d]^-\phi \\ \tilde D_{\phi}
 \ar@{|->}[r]^-{\chi_x} & \chi_x(\tilde D_\phi)=\phi(D_x)\,\,\,\,\,\,\,\,\,\,\,\,\,\,\,\,\,\,\,\,}$$

The composition:
$$G_a^\vee \overset{exp_{-D}\times f}\longrightarrow
{\mathbb Aut}\,\mathbb X\times \mathbb X \overset\circ \to \mathbb X$$
maps $\delta_\mu$ to zero, for all $\mu\in G_a^\vee$:
$$\begin{array}{rlcll} T_\mu G_a^\vee & \to
& T_{exp_{-D}(\mu)}{\mathbb Aut}\mathbb X\times T_{f(\mu)}\mathbb X & \to & T_{exp_{-D}(\mu)f(\mu)}\mathbb X\\
{\delta}_{\mu} & \mapsto & (-\tilde D_{exp_{-D}(\mu)},D_{f(\mu)}) & \mapsto & -exp_{-D}(\mu)D_{f(\mu)}+exp_{-D}(\mu)D_{f(\mu)}\end{array}$$ Since $\mathcal K=\{a\in \mathbb A_{G_a^\vee}\colon \delta_\mu(a)=0$, for all $\mu\in G_a^\vee\}$,
the composition $G_a^\vee
\overset{exp_{-D}\times f}\longrightarrow{\mathbb Aut}\mathbb X\times \mathbb X \to \mathbb X$,  is constant, that is, $exp_{-D}(\mu)f(\mu)=x$, and
$f(\mu)=exp_D(\mu)x=exp_{D,x}(\mu)$.

\end{proof}

\begin{theorem}[existence and uniqueness\,] \label{632} Let $X$ be a $K$-scheme, let $y \in X$ be a rational point and let $D$ be a vector field on $X$. Then,  $exp_{D,y}$ is
the only morphism $f\colon \Spec K[x]^* \to X$ such that $f(0)=y$ and $f({\delta}_{\mu})=D_{f(\mu)}$, for every point
$\mu\in \Spec K[x]^*(S)$.\end{theorem}

\begin{proof} By \ref{8} and \ref{63}, we only have to prove that given a morphism
$f\colon \Spec K[x]^* \to X$ such that $f(0)=y$ and $f({\delta}_{\mu})=D_{f(\mu)}$, for every point
$\mu\in G_a^\vee(S)\subset (\Spec K[x]^*)(S)$, then  $f({\delta}_{\mu})=D_{f(\mu)}$, for every $\mu\in\Spec K[x]^*(S)$. Consider the diagram
$$\xymatrix{A_{X,y} \ar[r]^-{f^*} \ar[d]^-D & K[x]^* \ar[d]^-\delta \\ A_{X,y} \ar[r]^-{f^*} & K[x]^*}$$
By the hypothesis $(f^*\circ D-\delta \circ f^*)(a)(\mu)=0$ for
all $a\in A_{X,y}$ and $\mu \in G_a^\vee(S)$. Then, $(f^*\circ
D-\delta \circ f^*)(a)=0$ and $f^*\circ D=\delta\circ f^*$. Hence,
$f({\delta}_{\mu})=D_{f(\mu)}$, for every $\mu\in\Spec K[x]^*(S)$.
\end{proof}

\subsection{Algebraic group associated with a vector field} ${}$

\medskip

Let $X$ be a $K$-scheme, let $D$ be a vector field on $X$ and let $\tilde D$ be the invariant field of $\mathbb End \mathbb X$ associated with the tangent vector $D$ at $\Id\in \mathbb End\, \mathbb X$.

\begin{proposition} \label{5.1} The analytic integral curve of $\tilde D$ passing through $\Id$ is $exp_D$, that is, $exp_{\tilde
D_{\Id}}= exp_D$.\end{proposition}

\begin{proof} $exp_D$ is a morphism of groups, so $exp_D({\delta}_{\mu})= {\tilde D}_{exp_D(\mu)}$. It follows from Theorem \ref{63} that
$exp_D=exp_{\tilde D_{\Id}}$.

\end{proof}

\begin{proposition} \label{5.2} The (closed) scheme-theoretic image of $exp_{D,x}$, $\Ima exp_{D,x}$, is the minimal subvariety of
$X$ tangent to $D$ at $x$.\end{proposition}

\begin{proof} Let $C$ be the minimal subvariety of $X$ tangent to $D$ at $x$. Let $$exp_{D,x}^C\colon \Spec K[x]^*\to C$$ be
the analytic integral curve associated with  $D_{|C}$ (passing through $x$). Obviously, the composition of
$exp_{D,x}^C$ with the inclusion $C\hookrightarrow X$ is $exp_{D,x}$, so $\Ima exp_{D,x}\subseteq C$. It is easy to check that
$\Ima exp_{D,x}$ is a subvariety  tangent to $D$ at $x$, so it is the minimal subvariety of $X$ tangent to $D$ at $x$.
\end{proof}

The ideal of functions vanishing at $\Ima exp_{D,x}$ is the greatest ideal whose germ at $x$ is the ideal of all
functions $f$ such that $D^nf(x)=0$ for every $n\geq 0$, if $\charc K=0$, or such that $D^{ip^n}f(x)=0$ for some $n\geq
0$ and every $0\leq i < p$, if $\charc K=p>0$.

\begin{corollary} Let $X$ be a projective variety. The (closed) scheme-theoretic image of $exp_D$ is the minimal closed subvariety of
${\mathbb Aut} X$ tangent to $\tilde D$ at $\Id$. $\Ima exp_D$ is a commutative group and the closure of $(\Ima exp_D)\cdot
x$ coincides with $\Ima exp_{D,x}$.\end{corollary}

\begin{proof} The scheme-theoretic image of $exp_D$ is the minimal closed subvariety of ${\mathbb Aut} X$ tangent to $\tilde D$ at
$\Id$, by \ref{5.1} and \ref{5.2}. $\Ima exp_D$ is a commutative group by \ref{1.3}. By definition, the composition
$$\xymatrix{ \Spec K[x]^* \ar[r]^-{exp_D} & {\mathbb Aut}\, X \ar[r] & X \\
& f \ar@{|->}[r] & f\cdot x = f(x) }$$ is $exp_{D,x}$. Therefore, $\Ima exp_{D,x}$ is the scheme-theoretic image of
$\Ima exp_D$, that is, the closure of $\Ima exp_D\cdot x$ coincides with $\Ima exp_{D,x}$.

\end{proof}

\begin{definition} \label{75} Let $X$ be a projective variety. We will say that $\Ima exp_D$ is the algebraic group associated with $D$.\end{definition}

The requirement of projectiveness is to assure the existence of the scheme ${\mathbb Aut}\,X$. It could be possible to define
the algebraic group $\Ima exp_D$ whenever $exp_D$ factors through a group scheme $G$ included in ${\mathbb Aut}\, X$ (it is
enough, for this, that $D\in T_{\Id} G$) because in that case we can also define $\Ima exp_D$.

\begin{theorem} Assume $\charc K=0$. Let $X=\Spec A$ be an affine
$K$-scheme. Then, $\Hom_{mon}(G_a,{\mathbb End}\,X)=\{D\in \Der
X\colon$ for each  $a\in A$  there exists an  $n\in\mathbb N$,
such that $D^n(a)=0\}$.
\end{theorem}

\begin{proof} Let $V$ be a vector space. Endowing $V$ with a $G_a$-module
structure is equivalent to endowing $\mathcal V$ with
a $\mathcal K[[x]]$-module structure. Endowing $\mathcal V$ with a 
$\mathcal K[[x]]$-module structure of  is equivalent to
defining an endomorphism $T\colon V\to V$ ($T=x\cdot$) such that
for each $v\in V$ there exists an $n\in\mathbb N$ so that
$T^n(v)=0$. Therefore,
$$\aligned   \left\{\begin{array}{l} T\in End_KV\colon \text{ for each } v\in V, \text{ there exists an }
n\in \mathbb N \\ \text{ such that } T^n(v)=0\end{array}\right\} &
=\Hom(G_a,{\mathbb End_KV})\\ T & \mapsto e^{x T}\endaligned$$
Arguing as in \ref{6.3}, we can prove this theorem.\end{proof}

\begin{theorem} Let $X=\Spec K[\xi_1,\ldots,\xi_n]$ be an affine
variety and assume $\charc K=0$. Let $D$ be a vector field of $X$.
Then, the morphism $exp_D\colon G_a^\vee\to {\mathbb Aut} X$
factors through an unipotent group (that is $G_a$) if and only if
there exists $m>>0$ such that $D^m(\xi_{i})=0$, for all $i$.

\end{theorem}

\begin{theorem} Assume $\charc K=p>0$ and let $X$ be a $K$-scheme.
It holds that
$$\Hom_{grp.} (\alpha_n^\vee,{\mathbb End}X)
=\{D\in \Der X\colon D^{p^n}=0\}$$
\end{theorem}

\begin{proof} 1. Let $V$ be a $K$-vector space. Endowing $V$ with  a $\alpha_n^\vee$-module
structure 
 is equivalent to endowing $V$ of a
 $K[x]/(x^{p^n})$-module structure. Then,

$$\{T\in \mathbb End_KV\colon T^{p^n}=0\}=\Hom_{funct}(\alpha_n^\vee,{\mathbb End_K V}), \, T\mapsto e^{xT}$$

2. Suppose, now, that $V$ is a $K$-algebra. Then, as in \ref{6.3},
$$\{D\in \Der_K(V,V)\colon D^{p^n}=0\}=\Hom_{grp}(\alpha_n^\vee, {\mathbb End_{K-alg} V})$$

3. Now, we can prove the theorem as in  \ref{55}.\end{proof}

\begin{theorem} Let $X$ be a complete $K$-variety and  $\charc K=p>0$.
Let $D$ be a vector field of $X$. Then, the morphism $exp_D\colon
G_a^\vee\to {\mathbb Aut} X$ factors through a unipotent group,
that is, $\alpha_n^\vee$  if and only if $D^{p^n}=0$.\end{theorem}

\subsection{Algebraic groups associated with the vector fields of $\mathbb P^n$} ${}$

\medskip

\label{76} Let $K$ be an algebraically closed field and consider the natural morphism $\pi\colon K^{n}\backslash 0\to
\mathbb P^{n-1}$. We are going to compute the algebraic group associated with a field on $\mathbb P ^{n-1}(K)$.

As it is well known, if $D$ is a vector field on $\mathbb P^{n-1}$, then $D=\pi D'$ for some vector field $D'=\sum_{ij} \lambda_{ij} x_i
\partial_{x_j}$ on $K^n\backslash 0$. We have the commutative diagram:

$$\xymatrix{ G_a^\vee \ar[rr]^-{exp_{D'}} \ar[rd] & & {\mathbb Aut}\,
K^n \\ & Gl_n(K) \ar@{^{(}->}[ru] & } \qquad \xymatrix{\mu \ar@{|->}[rr] \ar@{|->}[rd] & & e^{\mu D'} \\ & e^{\mu \cdot
(\lambda_{ij})} \ar@{|->}[ru] & }$$ The morphism $\Spec K[x]^* \to \Spec K[x_{ij}, \det(x_{ij})^{-1}]$, $\mu\mapsto
e^{\mu\cdot (\lambda_{ij})}$ is induced by $K[x_{ij}, \det(x_{ij})^{-1}]\to K[x]^*$, $x_{rs}\mapsto x_{rs}(e^{x\cdot
(\lambda_{ij})})$ (if $\charc K=0$ then $K[x]^*=K[[x]]$, if $\charc K=p>0$ then $K[x]^*=(K[[x_i]]/(x_i^p))_{i\in\mathbb N}$).

Therefore, the algebraic group associated with $D'$ is:
$$G=\Spec K[x_{rs}(e^{x\cdot (\lambda_{ij})}), e^{- x\cdot
\tr((\lambda_{ij}))}]$$

Changing the base of $K^n$, we can suppose that the matrix $(\lambda_{ij})$ is in its Jordan form, with eigenvalues
$\lambda_1,\ldots,\lambda_s$.

Now there are two cases, depending on the characteristic of $K$.

\begin{theorem} Let $K$ be an algebraically closed field of characteristic zero. Let $D=\pi(\sum_{ij}\lambda_{ij}x_i\partial_{x_j})$ be a vector field on $\mathbb P^{n-1}$ and let $G$ be its
associated algebraic group. Then,
$$G\simeq G_m^r\times G_a^\delta$$
where $r$ is the dimension of the $\mathbb Q$-affine space
generated by the eigenvalues of the matrix  $(\lambda_{ij})$ in
$\mathbb C$, $\delta=0$ in case the matrix is diagonalizable and
$\delta=1$ otherwise.\end{theorem}

\begin{proof} As the matrix $(\lambda_{ij})$ is in its Jordan form, it is easy to check that:
$$G=\Spec K[e^{x\lambda_1 },\ldots, e^{x\lambda_s }, e^{-x(\lambda_1+\cdots+\lambda_s)
}, \delta\cdot x]$$ where $\delta=0$ in case the Jordan matrix is diagonal  and $\delta=1$ otherwise. Moreover, if,
reordering, $\lambda_1,\ldots,\lambda_r$ is a base of the $\mathbb Z$-module generated by $\lambda_1,\ldots,\lambda_s$ in
$K$ (or, equivalently, $\lambda_1,\ldots,\lambda_r$ is a base of the $\mathbb Q$-vector space generated by
$\lambda_1,\ldots,\lambda_s$ in $K$) then, as $x,e^{x\lambda_1 },\ldots, e^{x\lambda_r }$ are algebraically
independent, we have that:
$$\aligned G&=\Spec K[e^{x\lambda_1 },\ldots, e^{x\lambda_s }, e^{-x(\lambda_1+\cdots+\lambda_s)
}, \delta\cdot x]\\ &=\Spec K[e^{x\lambda_1 },e^{-x\lambda_1 }]\otimes \cdots \otimes K[e^{x\lambda_r },e^{-x\lambda_r
}]\otimes K[\delta x] = G_m^r\times G_a^\delta\endaligned$$

\medskip But this group  $G$ is also the algebraic group associated with $D$. To see this, recall that two vector fields
$D'=\sum_{ij}\lambda_{ij}x_i\partial_{x_j}$ and $D'=\sum_{ij}\mu_{ij}x_i\partial_{x_j}$
on $K^n - \{ 0\}$ projet to the same vector field $D$ on $\mathbb{P}^{n-1}$ if and only if they differ on
$\lambda\cdot (\sum_i x_i \partial_{x_i})$, that is, if the matrix $(\mu_{ij})$ differs from
$(\lambda_{ij})$ in $\lambda\cdot\Id$. We can then assume that the matrix $M=(\lambda_{ij})$ has the eigenvalue zero,
and in this case the affine space generated by the eigenvalues is the vector space generated by the eigenvalues.

Let $v\in K^n$ be an eigenvector such that $M\cdot v=0$. If $H_v$
stands for the closed subgroup of $Gl_n(K)$ of the automorphisms
$h\in Gl_n(K)$ such that $h\cdot v=v$, then it is clear that the
morphism $exp_{D'}\colon \hat G_a \to Gl_n(K)$, $\mu\mapsto
e^{\mu\cdot M}$ factors through $H_v$. Analogously, let $H'_v$ be
the closed subgroup of the projectivities that let $\bar
v\in\mathbb P^{n-1}$ fixed.

At the Lie algebras, the natural morphism $Gl_n(K)\to
PGl_n(K)={\mathbb Aut}\, \mathbb P^{n-1}$ maps $D'$ to $D$.
Therefore, $exp_D$ coincides with the composition ${\hat
G_a}\overset{exp_{D'}}\to Gl_n(K)\to {\mathbb Aut}\, \mathbb
P^{n-1}$. As a consequence, we have:
$$\xymatrix{{\hat G_a} \ar[r]^-{exp_{D'}} \ar[rd] & Gl_n(K) \ar[r] &
PGl_n(K)\\ & H_v \ar@{=}[r] \ar@{^{(}->}[u] & H'_v
\ar@{^{(}->}[u]}$$ that shows that the algebraic group associated
with $D'$ is the same as that associated with $D$, so we are done.
\end{proof}

\begin{theorem} Let $K$ be an algebraically closed field of characteristic $p>0$. Let
$D=\pi(\sum_{ij}\lambda_{ij}x_i\partial_{x_j})$ be a vector field on $\mathbb P^{n-1}(K)$, and let $G$ be its associated
algebraic group. Then,
$$G\simeq \mu_1^r\times \alpha_{m+1}^\vee$$
where $r$ is the dimension of the $\mathbb Z/p\mathbb Z$-affine space generated by the eigenvalues of the matrix $(\lambda_{ij})$
in $K$, and $m$ is such that, if $s$ is the greatest of the orders of the Jordan boxes, then $p^m\leq s-1<p^{m+1}$ (if
$s=1$ we say that $m=-1$).\end{theorem}

\begin{proof} Everything is analogous to the previous theorem, except the calculation of the algebraic group $G$ associated with
$D' = \sum_{ij}\lambda_{ij}x_i{\partial_{x_j}} $ on $K^n - \{ 0\}$.

Reordering, let $\lambda_1,\ldots,\lambda_r$ be a base of the $\mathbb Z/p\mathbb Z$-vector space generated by
$\lambda_1,\ldots,\lambda_s$ in $K$. Let $s$ be the greatest of the orders of the Jordan boxes and $m\in\mathbb N$ such that
$p^m\leq s-1<p^{m+1}$ (if $s=1$ we say that $m=-1$). A similar computation to that used in the characteristic cero
case, shows that:
$$\aligned G&=\Spec K[e^{x\lambda_1 },\ldots, e^{x\lambda_s }, e^{-x(\lambda_1+\cdots+\lambda_s)
}, x_0,\ldots,x_{m}]/(x_0^p,\ldots,x_n^p,\ldots)\\ & = \mu_1^r\times \alpha_{m+1}^\vee\endaligned$$

\end{proof}

\section{Appendix} \label{Appendix}

\subsection{Tannakian Categories} ${}$

\medskip

\label{6.8} In this subsection we use Theorem \ref{Gmodulos} to derive the so called Tannaka's theorem (see \cite{Breen} and references therein for the standard treatment).

\medskip Let $K$ be a field.

\begin{definition} A \textit{neutralized $K$-linear category} $(\mathsf{C} , \omega)$ is an
abelian category $\mathsf{C}$ together with a ``fibre'' functor
$\omega\colon\mathsf{C}\rightsquigarrow {\rm Vect}_K$ into the
category of finite dimensional $K$-vector spaces such that $\omega$
is exact, additive and for every $X,X'\in \rm{Ob}(\mathsf C)$,
$$\Hom_{\mathsf{C}}(X,X')\subset \Hom_{K}(\omega(X),\omega(X'))$$
is a $K$-linear vector subspace.\end{definition}

A \textit{$K$-linear morphism} between neutralized $K$-linear categories $F
\colon ( \mathsf{C} , \omega ) \to (\bar{\mathsf{C}} ,
\bar{\omega})$ is an additive functor $F \colon \mathsf{C}
\rightsquigarrow \bar{\mathsf{C}}$ such that $\bar{\omega} \circ F
= \omega$.

 \begin{example} Let $A$ be finite a $K$-algebra. The category ${\rm Mod}_{A}$ of finitely generated modules over $A$ together with the forgetful functor is a neutralized $K$-linear category.

Recall also that morphisms of $K$-algebras $A \to B$ correspond to $K$-linear morphisms ${\rm Mod}_B \to {\rm Mod}_A$.
\end{example}

\medskip If  $(\mathsf{C} , \omega)$ is a neutralized $K$-linear category and $X
\in {\rm Ob}\, \mathsf{C}$ is an object, we will denote by $\langle
X \rangle$ the full subcategory of $\mathsf{C}$ whose objects are
(isomorphic to) quotients of subobjects of finite direct sums $X
\oplus \ldots \oplus X$.

By standard arguments, it can be proved the following:

\begin{theorem}[Main Theorem]  Let $\langle X \rangle$ be a neutralized $K$-linear category generated by an object $X$. There exists a (weak) equivalence of neutralized $K$-linear categories $\langle X \rangle \simeq {\rm Mod}_{A^X}$, where $A^X$ is a finite $K$-algebra unique up to isomorphisms.

Moreover, every $K$-linear morphism $F \colon \langle X \rangle \rightsquigarrow \langle \bar{X} \rangle$ induces a unique morphism of $K$-algebras $f \colon A^{\bar{X}} \to A^X $.\end{theorem}


A neutralized $K$-linear category $(\mathsf{C} , \omega)$ is said to {\sl{admit a set of generators}} if there exists a filtering set $I$ of objects in $\mathsf{C}$ such that:
$ \mathsf{C} = \underset{X \in I}{\underset{\to}\lim} \ \langle X \rangle .$

In this case, a standard argument passing to the limit allows to prove:
$$ \mathsf{C} = \lim_{\stackrel{\rightarrow}{X \in I}} \ \langle X \rangle \simeq \lim_{\rightarrow}\ {\rm Mod}_{A^X} =  {\rm Mod}_{{\mathcal C}^*}$$ where ${\mathcal C}^* := \lim \limits_{\leftarrow} \mathcal A^X$, that is a $\mathcal K$-algebra scheme and ${\rm Mod}_{{\mathcal C}^*}$ is the category of $\mathcal K$-coherent $\mathcal C^*$-modules.

Let $( \mathsf{C} , \omega )$ and $( \bar{\mathsf{C}} , \bar\omega )$ be two neutralized $K$-linear  categories that admit a set of generators.
Every $K$-linear morphism $F \colon ( \mathsf{C} , \omega ) \rightsquigarrow ( \bar{\mathsf{C}} , \bar{\omega} ) $ induces a unique morphism of $\mathcal K$-algebra schemes $f \colon \bar{{\mathcal C}}^* \to {\mathcal C}^* $.

\begin{definition} A tensor product on a neutralized $K$-linear category
$(\mathsf{C} , \omega)$ is a bilinear functor $\otimes \colon
\mathsf{C} \times \mathsf{C} \rightsquigarrow \mathsf{C}$ that
fits into the square:

$$\xymatrix{ {\mathsf C} \times {\mathsf C} \ar[d]^-{\omega\times \omega}
\ar[r]^-\otimes & {\mathsf C}\ar[d]^-\omega\\
{\rm Vect}_K \times {\rm Vect}_K \ar[r]^-{\otimes_K} & {\rm
Vect}_K}$$ (where the symbol $\otimes_K$ denotes the standard tensor
product on vector spaces) and satisfies:

a) Associativity and commutativity.

b) Unity. There exists an object $K$ together with functorial isomorphisms for every object $X$: $$X \otimes
{K} \simeq  X \simeq {K} \otimes X $$
that through $\omega$ become the natural identifications $\omega(X)\otimes_K K = \omega(X)=K\otimes \omega(X)$.

c) Duals. There exists a covariant additive functor
$\phantom{a}^\vee \colon \mathsf C \to  \mathsf{C}^\circ$,
satisfying:
$$\xymatrix{ {\mathsf C} \ar[r]^-{\phantom{a}^\vee} \ar[d]^-\omega &
{\mathsf C}^\circ \ar[d]^-{\omega^*}
\\ {\rm Vect}_K \ar[r]^-{*} & {\rm Vect}_K^\circ }$$ where $\omega^* (X)
:= \omega(X)^* $. There also exists functorial isomorphisms
$(X^\vee )^\vee = X$ and a morphism $ \mathcal{K} \to X \otimes
X^\vee $ such that via $\omega$ is the natural morphism $K\to
\omega(X)\otimes_K \omega(X)^*$.\end{definition}

\begin{definition} A {Tannakian category} neutralized  over $K$ is a
triple $(\mathsf{C} , \omega , \otimes)$ where $(\mathsf{C} ,
\omega )$ is a neutralized $K$-linear category that admits a set of generators and $\otimes$ is a tensor product
on $(\mathsf{C} , \omega )$.\end{definition}

Now it is not difficult to check that the existence of a tensor
product in a neutralized $K$-linear category $\mathsf{C} \simeq {\rm Mod}_{{\mathcal C}^*}$ amounts to the existence of a comultiplication on the algebra scheme ${\mathcal C}^*$. As a consequence:

\begin{theorem} \label{EsquemasAlgHopf} Let $(\mathsf{C} , \omega , \otimes)$ be a Tannakian category neutralized over $K$. There exists a unique (up to isomorphism)
cocommutative Hopf algebra $\mathcal K$-scheme ${\mathcal C}^*$ such that $(\mathsf{C} , \omega , \otimes)$ is equivalent to the
category ${\rm Mod}_{{\mathcal C}^*}$.
\end{theorem}

\begin{corollary}[Tannaka's Theorem] If $(\mathsf{C} , \omega , \otimes)$ is a
Tannakian category neutralized over $K$, then there exists a unique (up to isomorphism) affine
$K$-group scheme $G$ such that $(\mathsf{C} , \omega , \otimes)$ is
equivalent to the category of finite linear representations of $G$.
\end{corollary}

\begin{proof} By the previous theorem, there exists a scheme of Hopf algebras ${\mathcal C}^*$ such that $\mathsf{C} \simeq {\rm Mod}_{\mathcal C^*}$. If we define the affine group scheme $G := {\rm Spec}\, {C}$, then the statement follows from Theorem \ref{Gmodulos}.
\end{proof}

\subsection{Lie Algebras and Infinitesimal Formal Monoids in Characteristic Zero} ${}$

\medskip

\begin{notation} In this subsection all algebra schemes are assumed to be commutative and $R=K$ is assumed to be a field of characteristic zero.\end{notation}

Let $f\colon N\to M$ be a morphism of $R$-modules and let $f^*\colon \mathcal M^*\to \mathcal N^*$ be the dual morphism.   ${\rm Ker} f^*=\mathcal M_1^*$ and ${\rm Coker}\, f^*=\mathcal N_1^*$, where $M_1={\rm Coker}\, f$ and $N_1={\rm Ker}\, f$.

Let $\mathcal C^*$ be a commutative algebra scheme.   If $f^*\colon\mathcal M^*\to\mathcal N^*$ is a morphism of $\mathcal C^*$-modules then ${\rm Coker}\, f^*$ and ${\rm Ker}\, f^*$,  are $\mathcal C^*$-modules.

Let $\mathcal I_j^* \hookrightarrow \mathcal C^*$ ideal schemes and $m\colon {\mathcal I_1^*}\otimes \cdots\otimes \mathcal I_n^*\to\mathcal C^*$ the obvious multiplication morphism.
We denote by $\mathcal I_1^*\cdots \mathcal I_n^*=\mathcal J^*$ the module scheme closure  of ${\rm Im}\, m$  in $\mathcal C^*$, which is  an ideal scheme  of $\mathcal A^*$: the dual morphism of $m$, $c\colon C \to I_1\otimes\cdots\otimes  I_n$, is a morphism of
$\mathcal C^*$-modules and $J={\rm  Im}\, c$.

Given a functor of $\mathcal K$-modules ${\mathbb M}$ we will denote its $\mathcal K$-module scheme closure $\bar {\mathbb M}$.
 Observe that $\mathbb M^*(K)=\bar{\mathbb M}^*(K)$. Hence, $\bar {\mathbb M}=\mathcal N^*$, where $N=\mathbb M^*(K)$ (see \cite[2.7]{Amel}).
We say that a morphism of functors of $\mathcal K$-modules ${\mathbb M}\to {\mathbb N}$ is dense if $\bar {\mathbb M}\to \bar {\mathbb N}$ is surjective,  that is,
 if ${\mathbb N}^*(K)\to {\mathbb M}^*(K)$ is injective.

 We have
$${\mathcal I_1^*}\otimes \cdots\otimes \mathcal I_n^*\overset{\text{dense}}\to
\mathcal I_1^*\cdots \mathcal I_n^*\hookrightarrow \mathcal C^*$$
$\mathcal I_1^*\cdots \mathcal I_n^*=\mathcal I_1^*\cdot (\mathcal I_2^*\cdots \mathcal I_n^*)$: observe the diagram
$${\mathcal I_1^*}\otimes ({\mathcal I_2^*}\otimes\cdots\otimes \mathcal I_n^*)\overset{\text{dense}}\to {\mathcal I_1^*}\otimes \mathcal I_2^*\cdots \mathcal I_n^*\overset{\text{dense}}\to
\mathcal I_1^*\cdot(\mathcal I_2^*\cdots \mathcal I_n^*)\hookrightarrow \mathcal C^*$$
$\mathcal I_1^*\cdot (\mathcal I_2^*\cdot \mathcal I_3^*)=\mathcal I_1^*\cdot \mathcal I_2^*\cdot \mathcal I_3^*=(\mathcal I_1^*\cdot \mathcal I_2^*)\cdot \mathcal I_3^*$.

\begin{notation} Let $\mathbb M$ be an $\mathcal R$-module. We denote $S^n\mathbb M$ the functor of $\mathcal R$-modules
defined by $(S^n\mathbb M) (S):= S^n(\mathbb M(S))$ the $n$-th symmetric power of the $S$-module $\mathbb M(S)$. Let $S_n\mathbb M$ be the functor of modules defined by
$(S_n\mathbb M) (S):= (\mathbb M(S)\otimes_S\overset n\cdots\otimes_S\mathbb M(S))^{S_n}$. The natural morphism $S_n\mathbb M\to S^n\mathbb M$ is an isomorphism when $R$ is a $\mathbb Q$-algebra.\end{notation}

 Denote $\mathcal I^{*n}=\mathcal I^*\overset n\cdots\mathcal I^*$.
The composition $\mathcal I^*\otimes\overset n\cdots\otimes \mathcal I^*\to \mathcal I^{*n}\to \mathcal I^{*n}/\mathcal I^{*n+1}$ is dense and factors through
$S^n(\mathcal I^*/\mathcal I^{*2})$. Then, the morphism
$$\overline{{S}^n
({\mathcal I}^*/{\mathcal I}^{*2})}\overset m \to {\mathcal I}^{*n} / {\mathcal I}^{*n+1}$$
is  surjective.

Observe that $S^n\mathcal M$ is the quasi-coherent module associated to the $K$-module $S^nM$ and
$\overline{S^n\mathcal M^*}=(S^n\mathcal M)^*$, because
$$\aligned \Hom_{\mathcal K}(\overline{S^n{\mathcal M^*}},{\mathcal K}) & =\Hom_{\mathcal K}(S^n{\mathcal M^*},{\mathcal K})={\Hom}_{\mathcal K} (\mathcal M^* \otimes \stackrel{n}{\ldots} \otimes
\mathcal M^*, {\mathcal K} )^{S_n}\\ &=(M\otimes \stackrel{n}{\ldots}\otimes {M})^{S_n} = S^nM\endaligned$$

\begin{definition} Let ${\mathcal A}^*$ be a commutative bialgebra scheme, $e\colon \mathcal A^*\to \mathcal K$ its counit and $\mathcal I^*={\rm Ker}\, e$. We will say that $\mathbb G:={\rm Spec}\,\mathcal A^*$ is an infinitesimal formal monoid if $\mathcal A^*=\underset{\underset{i}{\longleftarrow}}{\lim}\,
\mathcal A^*/\mathcal I^{*i}$.
\end{definition}

\begin{theorem} \label{7.6} Let $\mathbb G={\rm Spec}\,\mathcal A^*$ be an infinitesimal formal monoid, $e\colon\mathcal A^*\to \mathcal K$ the unit of $\mathbb G$ and $\mathcal I^*={\rm Ker}\, e$. The natural morphism
$$\overline{S^n\mathcal (\mathcal I^*/\mathcal I^{*2})}\overset m\to \mathcal I^{*n}/\mathcal I^{*n+1}$$ is an isomorphism.
\end{theorem}

\begin{proof} Let us construct the inverse morphism  $\mathcal I^{*n}/\mathcal I^{*n+1}\to \overline{ S^n(\mathcal I^*/\mathcal I^{*2})}$ of $m$:
consider the multiplication morphism  $\mathbb G\times\overset n\cdots\times \mathbb G\to \mathbb G$,
$(g_1,\ldots,g_n)\mapsto g_1\cdots g_n$, which corresponds to the comultiplication morphism $c\colon \mathcal A^*\to \mathcal A^*\tilde\otimes\cdots \tilde\otimes \mathcal A^*$.

For any $f\in \mathcal I^*$ we have that $$c(f)=\sum_{j=1}^n 1\otimes \cdots \otimes\overset j f\otimes \cdots \otimes 1\,\,{\rm mod}\, \sum_{r\neq s}^n \mathcal A^*\tilde\otimes \cdots\tilde\otimes\underset r{\mathcal I}^*\tilde\otimes \cdots \tilde\otimes \underset s{\mathcal I}^*\tilde\otimes \cdots\tilde\otimes\mathcal A^*$$ because the classes of $c(f)$ and $\tilde f:=\sum_{j=1}^n 1\otimes \cdots \otimes\overset j f\otimes \cdots \otimes 1$  in $\mathcal A^*/\mathcal I^*\otimes \cdots \otimes\overset s {\mathcal A^*} \otimes \cdots \otimes \mathcal A^*/\mathcal I^*=\mathcal A^*$ are equal to $f$, for every $s$, so
$$c(f)-\tilde f\in \cap_{s=1}^n(\sum_{r\neq s}^n \mathcal A^*\tilde\otimes \cdots\tilde\otimes\underset r{\mathcal I}^*\tilde\otimes \cdots\tilde\otimes\mathcal A)=\sum_{r\neq s}^n \mathcal A^*\tilde\otimes \cdots\tilde\otimes\underset r{\mathcal I}^*\tilde\otimes \cdots \tilde\otimes \underset s{\mathcal I}^*\tilde\otimes \cdots\tilde\otimes\mathcal A^*$$
(for the latter equation recall $\mathcal A^*=\mathcal K\oplus \mathcal I^*$). Therefore, we obtain the morphism
$$\begin{array}{lll}  \mathcal I^{*n}/\mathcal I^{*n+1} & \overset{\bar c}\longrightarrow
&  \mathcal I^{*}/\mathcal I^{*2}\tilde\otimes \overset n\cdots \tilde\otimes \mathcal I^{*}/\mathcal I^{*2}\subset \mathcal A^{*}/\mathcal I^{*2}\tilde\otimes \overset n\cdots \tilde\otimes \mathcal A^{*}/\mathcal I^{*2}\\ f_1\cdots f_n & \mapsto & \overline{c(f_1\cdots f_n)}=\overline{c(f_1)\cdots c(f_n)}=\overline{\sum_{\sigma\in S_n}
f_{\sigma(1)}\otimes\cdots\otimes f_{\sigma(n)}}\end{array}$$
for every $f_1,\ldots,f_n\in \mathcal I^*$, that defines a morphism $\bar c\colon  \mathcal I^{*n}/\mathcal I^{*n+1}\to \overline{ S^n(\mathcal I^{*}/\mathcal I^{*2})}$. Now it can be checked that $\bar c\circ m\colon \overline{S^n(\mathcal I^*/\mathcal I^{*2})}\to \overline{ S^n(\mathcal I^*/\mathcal I^{*2})}$ is equal to the homothety  with scale factor $n!$.

\end{proof}

\begin{definition} If $A$ is a bialgebra, we say that an element is primitive if
$c(a)=a\otimes 1+1\otimes a$, where $c$ is the comultiplication of $A$.\end{definition}

It can be checked that  $a\in A$ is a primitive element if and only if $a\in T_e\mathbb G:={\rm Der}_{\mathcal K}(\mathcal A^*,\mathcal K)=\Hom_{\mathcal K}(\mathcal I^*/\mathcal I^{*2},\mathcal K)$ .

The inclusion $T_e\mathbb G\hookrightarrow A$ is a morphism of Lie algebras that extends to a morphism of algebras $U(T_e\mathbb G)\to A$, where $U(T_e\mathbb G)$ is the universal algebra of $T_e\mathbb G$.

Let $ L$ be a Lie algebra. $U(L)$ is a quotient of the tensorial algebra of $L$, $T^\cdot L$. It is easy to see, (\cite[I.III.4.]{Serre}) that $S^\cdot L$ has a surjective morphism onto the graduated algebra by the filtration of $U(L)$, $\{U(L)_n:=\left[\oplus_{i\leq n}T^iL\right]\}$.

Let $\mathbb G={\rm Spec}\, \mathcal A^*$  be an infinitesimal formal group. Let us denote $\mathcal A_n=(\mathcal A^*/\mathcal I^{*n+1})^*$.
The equality $\mathcal A^*=\underset{\underset{i}{\longleftarrow}}{\lim}\,
\mathcal A^*/\mathcal I^{*i}$ is equivalent to the equality $\mathcal A =\underset{\underset{i}{\longrightarrow}}{\lim}\,
\mathcal A_i$. Observe that $A_i\cdot A_j\subseteq A_{i+j}$:
let $c\colon \mathcal A^*\to \mathcal A^*\tilde\otimes \mathcal A^*$ be the comultiplication. Then, $c({\mathcal I^*})\subset \mathcal I^*\tilde\otimes \mathcal A^*+\mathcal A^*\tilde\otimes \mathcal I^*$, so that $c({\mathcal I}^{*i+j+1})\subseteq {\mathcal I}^{*i+1}\tilde\otimes{\mathcal A^*+\mathcal A^*\tilde\otimes \mathcal I}^{*j+1}$. The dual morphism of
$$\mathcal A^*/\mathcal I^{*i+j+1}\overset c\to \mathcal A^*/\mathcal I^{*i+1}\tilde\otimes
\mathcal A^*/\mathcal I^{*j+1}
$$ is the multiplication morphism $A_i\otimes A_j\to A_{i+j}$. The morphism $U(L)\to A$ maps $U(L)_1$ into  $A_1$, so $U(L)_n$ maps into $A_n$. Lastly, it is easy to check that  $\mathcal A_n/\mathcal A_{n-1}=(\mathcal I^{*n}/\mathcal I^{*n+1})^*=:\mathcal L_n$.

\begin{theorem} \label{6.12} Let $\mathbb G={\rm Spec}\, \mathcal A^*$ an infinitesimal formal group, and write $L:=T_e\mathbb G$. Then,
\begin{enumerate}

\item $U(L)\to A$ is an isomorphism of bialgebras.

\item  The morphism $U(L)_n/U(L)_{n-1}\to A_n/A_{n-1}$ is an isomorphism.

\item  $L$ is the module of primitive elements of $U(L)$ and $U(L)_n/U(L)_{n-1}=S^nL$.
 \end{enumerate}\end{theorem}

\begin{proof}  From the commutative diagram
$$\xymatrix{S^nL \ar@{=}[r]^-{\text{\ref{7.6}}} \ar[rd]_-{\text{surj}} & L_n \\ & U(L)_n/U(L)_{n-1} \ar[u]}$$
it easily follows $(2)$. By induction on $n$, it is easy to see that $U(L)_n\to A_n$ is an isomorphism, then $U(L)\to A$ is an isomorphism.

Moreover, $U(L)\to A$ is a morphism of coalgebras because it maps $L$, that are primitive elements of $U(L)$, into primitive elements of $A$ and $U(L)$ is generated algebraically by $L$. Finally, the module of primitive elements of $A$ is $L$, so the module of primitive elements of  $U(L)$ is precisely $L$.

\end{proof}

\begin{corollary} Let $\mathbb G={\rm Spec}\, \mathcal A^*$, ${\mathbb G'}={\rm Spec}\,\mathcal B^*$ be infinitesimal formal groups. Then,
$${\rm Hom}_{grp}(\mathbb G,{\mathbb G'})={\rm Hom}_{Lie}(T_e\mathbb G,T_e{\mathbb G'})$$\end{corollary}

\begin{proof} It follows from:

$$\aligned {\rm Hom}_{grp}(\mathbb G,{\mathbb G'}) & ={\rm Hom}_{bialg}({\mathcal B^*},{\mathcal A^*})={\rm Hom}_{bialg}
(A,B)\\ & ={\rm Hom}_{bialg}(U(T_e\mathbb G),U(T_e{\mathbb G'}))={\rm Hom}_{Lie}(T_e\mathbb G,T_e{\mathbb G'})\endaligned$$

\end{proof}

\begin{note} If $L$ is a Lie algebra, consider $\mathbb G=\Spec \mathcal U(\mathcal L)^*$. Let $\bar L=T_e\mathbb G$, that is, the primitive elements of  $U(L)$. We have a natural morphism $L\to \bar L$. With adequate basis in $L$ and $\bar L$ we have the commutative diagram:
$$\xymatrix{ S^\cdot L \ar[d]^-{surj} \ar[r] & S^\cdot \bar L \ar@{=}[d]^-{\text{\ref{6.12}(2)}} \\ U(L)\ar@{=}[r]^-{\text{\ref{6.12}(1)}} & U(\bar L)}$$ that allows to prove that the morphism  $L\to \bar L$ is surjective.

Let us also outline very briefly that the morphism  $L\to \bar L$, $D\mapsto\bar D$ is injective (see \cite[5.4]{Serre}). We only have  to prove that there exists a faithful linear representation of $L$, since
$U(L)=U(\bar L)$. If $L$ is commutative, then $S^\cdot L=U(L)$ and the morphism $L\to U(L)$ is injective. Let $Z$ be the kernel of the surjection $L\to \bar L$ (notice that $[L,Z]=0$).
 Let $\mathbb G_Z$ and $\mathbb G_{\bar L}$ be the formal groups associated to  $Z$ and $\bar L$. It is enough to see that there exists a morphism of Lie algebras $L\to {\rm Der}_K(\mathbb G_Z\times \mathbb G_{\bar L})$ injective. To do that, it is enough to prove that there exists a section of Lie algebras $w\colon \bar L\otimes_K U(\bar L)^*\to L\otimes_K U(\bar L)^*$ of the natural surjection $L\otimes_K U(\bar L)^*\to \bar L\otimes_K U(\bar L)^*$. Let $s\colon \bar L\to L$ be any $K$-linear section. It can be checked that the 2-form of $\mathbb G_{\bar L}$ with values in $Z$, $w_2\colon \bar L\times \bar L\to Z$, $w_2(\bar D,\bar D')=s([\bar D,\bar D'])-[D,D']$ is closed. By  Poincare Lemma, there exists a 1-form of $\mathbb G_{\bar L}$ with values in $Z$, $w'\colon \bar L\otimes_K U(\bar L)^*\to Z\otimes_K U(\bar L)^*$, such that  $dw'=w_2$. The section of Lie algebras that we were looking for is $w=s+w'$.
\end{note}

\begin{theorem} The category of infinitesimal formal groups is equivalent to the category of Lie algebras.\end{theorem}

\begin{proof} The functors giving the equivalence assign to each infinitesimal formal group
$\mathbb G$ its tangent space at the identity $T_e\mathbb G$ and to each Lie algebra $L$, the group ${\rm Spec}\,\mathcal U(L)^*$.\end{proof}

\begin{corollary} The category of linear representations of an infinitesimal formal group
$\mathbb G$ is equivalent to the category of linear representations of its Lie algebra $T_e\mathbb G$.\end{corollary}

\begin{proof} The category of linear representations of the formal group
$\mathbb G={\rm Spec}\, \mathcal A^*$ is equivalent to the category of $A$-modules, that is equivalent to the category of linear representations of the Lie algebra $T_e\mathbb G$, because $A$ is the universal algebra associated to $T_e\mathbb G$.\end{proof}

Let $G={\rm Spec}\, A$ be an affine $K$-group scheme and $I_e$ the ideal of functions that vanish at the identity element of $G$. Let
$J$ be the set of ideals of finite codimension of $A$ that are included in $I_e$ and let us denote ${\rm Dist}(G) :=\underset{{I\in J}}{\underset\rightarrow\lim} (A/I)^*$.

\begin{corollary} Let $G={\rm Spec}\, A$ be an affine $K$-group scheme. There exists a canonical isomorphism of bialgebras: $$U(T_eG)={\rm Dist}\, G$$ Therefore, $U(T_eG)^*=\hat A$ and the infinitesimal formal group associated to $T_eG$ is $\hat G$.\end{corollary}

\begin{proof} Let $\hat{\mathcal A}:=\underset{I\in J}{\underset\leftarrow\lim} \mathcal A/\mathcal I$
and $\hat G={\rm Spec}\, \hat{\mathcal A}$. Observe that ${\rm Hom}_{\mathcal K}(\hat{\mathcal A}, \mathcal K)={\rm Dist}\, G$. Moreover,
$$T_eG={\rm Hom}_{{\rm Spec}\, K}({\rm Spec}\, K[x]/(x^2), G)={\rm Hom}_{{\rm Spec}\, \mathcal K}({\rm Spec}\, \mathcal K[x]/(x^2), \hat G)=T_e\hat G$$
Therefore, by Theorem
\ref{6.12}, ${\rm Dist}\, G=U(T_e\hat G)=U(T_eG)$. \end{proof}

(See \cite[III.6.1]{Demazure}, where $G$ is algebraic).

\begin{corollary} If $G={\rm Spec}\, A$ is a commutative unipotent $K$-group, then it is isomorphic to $\mathcal V^*$, where $V=T_eG^\vee$.\end{corollary}

\begin{proof} $G$ is a commutative unipotent $K$-group if and only if $G^\vee$ is a commutative infinitesimal formal group. By Theorem \ref{6.12}, $G^\vee={\rm Spec}\, \mathcal (U(T_eG^\vee))^*$. As $T_eG^\vee\subset A$ is a trivial Lie algebra,
$G={\rm Spec}\, U(T_eG^\vee)= {\rm Spec}\, S^\cdot (T_eG^\vee)=\mathcal V^*$.\end{proof}

\subsection{Another examples of affine functors} ${}$

\medskip

\begin{definition} Let $\mathbb M$ be a functor of $\mathcal R$-modules and
let $\mathbb A_{\mathbb M}$ be the functor of functions of $\mathbb X=\mathbb M$. We define
$$\left[\mathbb A_{\mathbb M}\right]_n:=\{F\in \mathbb A_{\mathbb M}\colon F(\lambda\cdot m)=
\lambda^n\cdot F(m)\text{, for all }\lambda\in\mathcal R\text{ and }m\in\mathbb M\}.$$
\end{definition}

\begin{proposition} Let $\mathbb M$ be a functor of $\mathcal R$-modules. Then,
$$\oplus_n \left[\mathbb A_{\mathbb M}\right]_n \subseteq \mathbb A_{\mathbb M} \subseteq \prod_n \left[\mathbb A_{\mathbb M}\right]_n,$$ which are inclusions of $\prod_{\mathbb N}\mathcal R$-modules. Specifically,
$$\aligned \mathbb A_{\mathbb M}=\{\sum_n F_n\in \prod_n \left[\mathbb A_{\mathbb M}\right]_n\colon &
\text{for each $m\in\mathbb M$ there exists $r\in\mathbb N$ such that } \\ & F_n(m)=0,
\text{ for all } n>r\}\endaligned$$
\end{proposition}

\begin{proof} Given $m\in\mathbb M$, let $G^m\colon \mathcal R\to \mathcal R$ be defined
by $G^m(\lambda)=F(\lambda m)$. Then, $G^m(x)=\sum_n r^m_n x^n\in \mathbb A_{\mathcal R}=\mathcal R[x]$ ($r^m_n=0$, for all $n>>0$) and $F(\lambda m)=\sum_n r^m_n \lambda^n$. Let
$F_n\colon \mathbb M\to \mathcal R$ be defined by $F_n(m)=r^m_n$ (given $m$, $F_n(m)=0$,  for all $n>>0$).
Observe that $F(\lambda (\mu m))=\sum_n r^m_n (\lambda \mu)^n=\sum_n (r^m_n\mu^n)\cdot \lambda^n$, then $F_n(\mu m)=r^m_n\mu^n=\mu^n F_n(m)$ and $F_n\in \left[\mathbb A_{\mathbb M}\right]_n$. Moreover, $F(m)=\sum_n r^m_n=\sum_n F_n(m)$.

Finally, let $(w_n)\in  \prod_n \left[\mathbb A_{\mathbb M}\right]_n$ such that for each $m\in\mathbb M$ there exists $r\in\mathbb N$ so that   $w_n(m)=0,$ for all  $n>r$. If $F:=\sum_n w_n=0$ then $w_n=0$, for all $n$: $0=(\sum_n w_n)(\lambda m)=
\sum_n w_n(m)\lambda^n$, for all $\lambda$, then $w_n(m)=0$ for all $n$, and $w_n=0$
for all $n$.

\end{proof}

\begin{proposition} \label{uff} Let $R$ be a commutative $\mathbb Q$-algebra. Let $\mathbb M$ be a functor of $\mathcal R$-modules. Then,

$$\left[\mathbb A_{\mathbb M}\right]_n=(S_n\mathbb M)^*.$$

\end{proposition}

\begin{proof} 

Let $i\colon \mathbb M\to S_n\mathbb M$ be the morphism of
functors of sets defined by $i(m)=m\otimes\overset n\cdots \otimes
m$ and $i^*\colon \mathbb A_{S_n \mathbb M}\to \mathbb A_{\mathbb
M}$ the morphism induced over the functors of functions. If $w\in
(S_n\mathbb M)^*$ then $i^*(w)=w\circ i\in \left[\mathbb
A_{\mathbb M}\right]_n$. If $i^*(w)=w\circ i=0$ then $w=0$: By
hypothesis, $w(m\otimes\overset n\cdots\otimes m)=0$, for all
$m\in\mathbb M$. Given $m_1,\ldots,m_n\in\mathbb M$ and $m=\sum_i
\lambda_i m_i$ ($\lambda_i$  being variables), then
$0=w(m\otimes\overset n\cdots \otimes m)=\sum_{|\alpha|=n}
\lambda^\alpha w(m^\alpha)$. Hence, $w(m^\alpha)=0$, and $w=0$.

Now,
let $F\in\left[ \mathbb A_{\mathbb M}\right]_n$. Given $m_1,\ldots,m_r\in\mathbb M$
let $G\colon \mathcal R^r\to \mathcal R$ be defined by $G((\lambda_i))=F(\sum_i\lambda_im_i)$. Then, $G=\sum_{|\alpha|=n}
r_\alpha x^\alpha\in \left[\mathbb A_{\mathcal R^r}\right]_n=\left[\mathcal R[x_1,\ldots,x_r]\right]_n$. That is,
$F(\sum_i\lambda_im_i)=\sum_{|\alpha|=n}
r_\alpha \lambda^\alpha$. Let $$F_n\colon \mathbb M\times\overset n\cdots \times\mathbb M\to \mathcal R,\, F_n(m_1,\ldots,m_n):=r_{(1,\overset n\ldots,1)}$$

Let us check that $F_n$ is a symmetric $n$-multilinear  mapping of $\mathbb M^n$, that is, $F_n\in  (S^n\mathbb M)^*$. Obviously $F_n$ is symmetric.
$F(\lambda_1 m_1+\cdots+\lambda_n\mu m_n)=\sum_{|\alpha|=n} r_\alpha \mu^{\alpha_n} \lambda^\alpha$, then $F_n(m_1,\cdots,\mu m_n)=\mu\cdot F_n(m_1,\ldots,m_n)$.
Let us write
$$\aligned F(\lambda_1m_1+\cdots+\lambda_{n+1}m_{n+1})& =\sum_{|\beta|=n} a_\beta \lambda^\beta\\ & =
\lambda_1\cdots\lambda_{n-1}\cdot (a_{(1,\ldots,1,1,0)}\lambda_n+a_{(1,\ldots,1,0,1)}\lambda_{n+1})\\ & +
\sum_{|\beta|=n, (\beta_1,\ldots,\beta_{n-1})\neq (1,\ldots,1)} a_\beta \lambda^\beta\endaligned$$
Considering $\lambda_{n+1}=0$ we obtain $F_n(m_1,\ldots,m_{n-1},m_n)=a_{(1,\ldots,1,1,0)}$.
Considering $\lambda_{n}=0$ we obtain $F_n(m_1,\ldots,m_{n-1},m_{n+1})=a_{(1,\ldots,1,0,1)}$.
Considering $\lambda_{n}=\lambda_{n+1}$ we obtain $F_n(m_1,\ldots,m_n+m_{n+1})=a_{(1,\ldots,1,1,0)}+a_{(1,\ldots,1,0,1)}$.
Hence $F_n$ is linear.

Let $\tilde F_n\colon  \mathbb M\otimes\overset n\cdots
\otimes\mathbb M\to \mathcal R$ be the morphism defined by $F_n$.
 Let us prove that $n!\cdot F$ is the composite morphism
$$\mathbb M\overset i\to S_n\mathbb M\subset \mathbb M\otimes\overset n\cdots
\otimes\mathbb M\overset{\tilde F_n}\to \mathcal R$$ Given
$m\in\mathbb M$, write $a_n:=F(m)$. Then, $F(\lambda\cdot
m)=\lambda^na_n$.
$F(\lambda_1m+\cdots+\lambda_nm)=F((\lambda_1+\cdots+\lambda_n)m)=
(\lambda_1+\cdots+\lambda_n)^n
a_n=n!(\lambda_1\cdots\lambda_n)a_n+\cdots$, hence, $(\tilde
F_n\circ i)(m)=F_n(m,\overset n\cdots, m)=n!a_n=n!\cdot F(m)$.

%
%
%
%

\end{proof}

\begin{proposition} \label{uf2} Let $M$ be a flat $R$-module. Then,
$$\left[\mathbb A_{\mathcal M}\right]_n=(S_n\mathcal M)^*,$$
which is a reflexive functor of modules.
\end{proposition}

\begin{proof}  By Govorov-Lazard Theorem (\cite[A6.6]{eisenbud}), $M$ is a direct limit of free modules of finite type, $M=\ilim{i} V_i$. Then,
$$\left[\mathbb A_{\mathcal M}\right]_n=\plim{i} \left[\mathbb A_{\mathcal V_i}\right]_n=
\plim{i}((S_n\mathcal V_i)^*)=(\ilim{i} S_n\mathcal V_i)^*=(S_n\mathcal M)^*$$
Observe that $S_n\mathcal M$ is a quasi-coherent module because it is a direct
limit of  quasi-coherent modules. Hence, $(S_n\mathcal M)^*$ is reflexive.

\end{proof}

\begin{proposition} \label{7.66} Let $R$ be a commutative $\mathbb Q$-algebra.
 Let $\mathcal M$  be a quasi-coherent $\mathcal R$-module, then $\mathbb X=\mathcal M$ is an affine functor.\end{proposition}

\begin{proof} By \cite[5.1]{navarro} and Lemma \ref{uff}, $\mathbb A_{\mathcal M}$ is reflexive.

 A morphism $\phi\colon \mathbb A_{\mathcal M}\to \mathcal R$ is determined by the restriction of $\phi$ on $\oplus_n (S_n\mathcal M)^*$ (see \cite[5.1]{navarro}). Observe that $(S^\cdot \mathcal M^*)^{**}=(\prod_n S_n\mathcal M)^*=\oplus_n (S_n\mathcal M)^*$. Then,

$$\aligned \Spec \mathbb A_{\mathcal M} & \subseteq \mathbb Hom_{\mathcal R-alg}(\oplus_n (S_n\mathcal M)^*,\mathcal R)\overset{\text{\ref{n2.16}}}=
\mathbb Hom_{\mathcal R-alg}(S^\cdot \mathcal M^*,\mathcal R)=\mathbb Hom_{\mathcal R}({\mathcal M}^*,\mathcal R)\\ & =\mathcal M\endaligned$$
The composition of the  natural morphism $\mathcal M\to \Spec \mathbb A_{\mathcal M}$ with the inclusion $ \Spec \mathbb A_{\mathcal M} \subset \mathcal M$ is the identity morphism.  Therefore, $\Spec \mathbb A_{\mathcal M}=\mathcal M$.
\end{proof}

\begin{proposition} \label{uf22} Let $M$ be a flat $R$-module. Then,
$\mathcal M$ is an affine functor.

\end{proposition}

\begin{proof}
Proceed as in the proof of Proposition \ref{7.66}.
\end{proof}

\begin{theorem} \label{pro314}Let $\mathbb M\in\mathfrak F$. Then, $\mathbb M$ is an affine functor.\end{theorem}

\begin{proof} There exist inclusions $\oplus_I\mathcal R\subseteq \mathbb M\subseteq \prod_I\mathcal R$, which are morphisms of $\prod_I\mathcal R$-modules.

Let us prove that $\mathbb A_{\mathbb M}\subseteq \mathbb A_{\oplus_I\mathcal R}$: Let $f\in \mathbb A_{\mathbb M}$ such that $f_{|\oplus_I\mathcal R}=0$. Given $m=(r_i)_{i\in I}\in\mathbb M\subseteq \prod_I\mathcal R$ let
$F\colon \prod_I\mathcal R\to\mathcal R$ be defined by
$F((\lambda_i))=f((\lambda_ir_i))$. Then, $F\in \mathbb A_{\prod_I\mathcal R}=
S^\cdot (\oplus_I \mathcal R)=\mathcal R[x_i]_{i\in I}$. Hence, there exists a finite subset $J\subseteq I$ such that $F\in \mathcal R[x_j]_{j\in J}$. Therefore
$F((\lambda_ir_i)_{i\in I})=F((\lambda_jr_j)_{j\in J})$ and
$f((r_i)_{i\in I})=f((r_j)_{j\in J})=0$, that is, $f=0$.

In conclusion we have
$$\underset{n\in\mathbb N}\oplus S^n(\oplus_I\mathcal R)=\mathbb A_{\prod_I \mathcal R}\subseteq\mathbb A_{\mathbb M}\subseteq \mathbb A_{\oplus_I\mathcal R}\overset{\text{\ref{uf2}}}\subseteq
\prod_{n\in\mathbb N} (S_n(\oplus_I\mathcal R))^*$$

Then, $\mathbb A_{\mathbb M}\in \mathfrak F$ and it is reflexive.

Any morphism of $\mathcal R$-modules $\phi\colon \mathbb A_{\mathbb M}\to\mathcal R$
is determined by its restriction to $S^\cdot(\oplus_I\mathcal R)$, by \cite[5.1]{navarro}.
Then, any morphism of $\mathcal R$-algebras $\varphi\colon \mathbb A_{\mathbb M}\to\mathcal R$ is determined by its restriction to $\oplus_I\mathcal R$. Since
$\oplus_I\mathcal R=(\prod_I\mathcal R)^*\subseteq \mathbb M^*$, $\varphi$ is determined
by its restriction to $\mathbb M^*$. Hence, $\Spec\mathbb A_{\mathbb M}\subseteq \mathbb M^{**}=\mathbb M$. Then, $\Spec\mathbb A_{\mathbb M}=\mathbb M$.

\end{proof}

\end{document}